\documentclass[10pt]{amsart}
\usepackage{amssymb}
\usepackage[dvips]{graphicx}
\usepackage{psfrag}

\theoremstyle{plain}
\newtheorem*{theorem A}{Theorem A}
\newtheorem*{theorem B}{Theorem B}
\newtheorem*{theorem C}{Theorem C}
\newtheorem*{theorem C'}{Theorem C'}
\newtheorem*{main lemma}{Main Lemma}
\newtheorem{theorem}{Theorem}[section]
\newtheorem{proposition}{Proposition}[section]
\newtheorem{lemma}{Lemma}[section]

\theoremstyle{definition}
\newtheorem*{definition}{Definition}

\theoremstyle{remark}
\newtheorem*{remark}{Remark}
\newtheorem*{claim}{Claim}
\newtheorem*{claim proof}{Proof of the claim}

\def\Diff{\mathrm{Diff}_{\mu}^{1}(M)}
\def\sldoisr{\mathrm{SL}(2,\mathbb{R})}
\def\RR{\mathbb{R}}
\def\ZZ{\mathbb{Z}}
\def\NN{\mathbb{N}}
\def\BB{\mathcal{B}}
\def\UU{\mathcal{U}}
\def\eps{\varepsilon}
\def\ang{\measuredangle}

\begin{document}

%%%%%%%%%%%%%%%%%%%%%%%%%%%%%%%%%%%%%%%%%%%%%%%
\title{Genericity of zero Lyapunov exponents}
\author{Jairo Bochi}
\date{September 15, 2000. Revised on \today}
\thanks{Work supported by FAPERJ and CNPq.} 
\address{IMPA -- Estr. D. Castorina 110 \newline
22460-320 Rio de Janeiro -- Brazil.}
\email{bochi@impa.br}
%\subjclass[2001]{12}

\begin{abstract}
We show that, for any compact surface, there is a residual (dense $G_{\delta}$)
set of $C^{1}$ area preserving diffeomorphisms which either are Anosov
or have zero Lyapunov exponents a.e. This result was announced by R.
Ma\~{n}\'{e}, but no proof was available. We also show that for any fixed
ergodic dynamical system over a compact space, there is a residual set of
continuous $\sldoisr$-cocycles which either are uniformly
hyperbolic or have zero exponents a.e.
\end{abstract}

\maketitle

%%%%%%%%%%%%%%%%%%%%%%%%%%%%%%%%%%%%%
%%%%%%%%%%%%%%%%%%%%%%%%%%%%%%%%%%%%%
\section{Introduction}
%%%%%%%%%%%%%%%%%%%%%%%%%%%%%%%%%%%%%
%%%%%%%%%%%%%%%%%%%%%%%%%%%%%%%%%%%%%

Let $M$ be a compact connected Riemannian two-dimensional $C^{\infty }$
manifold without boundary and let $\mu$ be its normalized area. Denote by 
$\Diff$ the set of all $\mu $-preserving $C^{1}$
diffeomorphisms endowed with the $C^{1}$ topology.

For $f\in \Diff$, the number 
\[
\lambda ^{+}(f,x)=\lim_{n\rightarrow +\infty }\frac{1}{n}\log \left\|
Df_{x}^{n}\right\| , 
\]
called the upper Lyapunov exponent, exists for almost every $x\in M$. Besides, $%
\lambda ^{+}(f,x)\geq 0$. Our main result is the following:

\begin{theorem A}
There exists a residual subset 
$\mathcal{R}\subset \Diff$ such that every
$f\in \mathcal{R}$ either is Anosov or $\lambda ^{+}(f,x)=0$
for $\mu $-almost every $x$.
\end{theorem A}

This theorem was announced by Ricardo Ma\~{n}\'{e} (1948-1995) around 1983.
He also announced its generalization to symplectic manifolds (with a
somewhat more elaborate statement) in \cite{Mane2}. Although his proofs
have never been published, a sketch of a proof of Theorem~A appeared in
1995, \cite{Mane1}. We exploited ideas outlined there, together with new
ingredients, to prove the theorem in the present paper.

%Roughly speaking, the proof goes as follows: 
%In the absence of uniform hyperbolicity, the expanding and contracting
%directions can be mixed. Doing this we construct local (i.e. along orbit segments)
%perturbations with small ``finite-time-exponent''.
%To extend the construction to the whole surface, a recurrence analysis is needed.

We recall that the only surface that admits Anosov
diffeomorphisms is the torus, see \cite{Franks}.

It is interesting to compare our results with another $C^{1}$-generic
dichotomy for area preserving diffeomorphisms that was obtained by Newhouse 
\cite{Newhouse}: A generic diffeomorphism either is Anosov or the set of
elliptic periodic points is dense in the surface.

We will indicate by $\mathrm{LE}(f)$ the ``integrated Lyapunov exponent'' of 
$f\in \Diff$, that is, 
\[
\mathrm{LE}(f)=\int_{M}\lambda ^{+}(f)d\mu . 
\]
Recall \cite{DesigRuelle} Ruelle's inequality $h_{\mu }(f)\leq \mathrm{LE}(f),$ where $h_{\mu
}(f)$ is the metric entropy of $f$.
As a corollary of the proof of theorem A, we obtain that the functions 
\[
\mathrm{LE}\mbox{ and }h_{\mu }:\Diff\rightarrow
\lbrack 0,\infty ) 
\]
are not continuous.

A more general setting to study Lyapunov exponents consists on
linear cocycles. In section~5, we prove the analogue of
theorem~A for this context. More precisely, let $X$ be a compact space, let $%
T:X\rightarrow X$ be a homeomorphism and let $\mu $ be some ergodic measure
for $T$. Given a continuous matrix function (called a linear cocycle) $%
A:X\rightarrow \sldoisr$, the following limit exists:
\[
\mathrm{LE}(A)=\lim_{n\rightarrow +\infty }\frac{1}{n}\log \left\|
A(T^{n-1}x)\cdots A(x)\right\| \geq 0
\]
for $\mu $-a.e. $x\in X$. We prove the following:

\begin{theorem C}
If $T$ is ergodic then there
is a residual set $\mathcal{R}\subset C^0(X,\sldoisr)$
such that for every $A\in \mathcal{R}$, either $A$ is
uniformly hyperbolic or $\mathrm{LE}(A)=0$.
\end{theorem C}

It is natural to ask whether theorem A extends, for instance, to the $C^{2}$
topology or whether theorem C extends to the $C^{1}$ topology. The answer to
the first question is unknown, but the answer to the second is negative. For
instance, Young \cite{LSYoung} exhibits open subsets of $C^{1}(X,\sldoisr)$
consisting of nonuniformly hyperbolic cocycles with positive
exponent, where the base transformations are linear automorphisms of the
two-torus.

Deciding whether a diffeomorphism (or a cocycle) has non-zero exponents is,
in general, very hard. Our results above provide some explanation for that
fact: positive Lyapunov exponents is \textit{not} a $C^{1}$-open condition.

A related very natural question in this setting is whether systems with
non-zero Lyapunov exponents are typical, in a measure-theoretic sense. At
this point, there is no general theorem in this direction.

%%%%%%%%%%%%%%%%%%%%%%%%%%%%%%%%%%%%%%%%%%%%%%%
%%%%%%%%%%%%%%%%%%%%%%%%%%%%%%%%%%%%%%%%%%%%%%%
\section{Preliminaries}
%%%%%%%%%%%%%%%%%%%%%%%%%%%%%%%%%%%%%%%%%%%%%%%
%%%%%%%%%%%%%%%%%%%%%%%%%%%%%%%%%%%%%%%%%%%%%%%

The manifold $M$ and the measure $\mu$ will be fixed from here until
the end of section~\ref{s.prova teor A}.
Also, ``a.e.'' will mean ``$\mu$-almost every''.

%%%%%%%%%%%%%%%%%%%%%%%%%%%%%%%%%%%%%%%%%%%%%%%%%%%%%%%%%%%
\subsection{Oseledets' theorem and Lyapunov Exponents}
%%%%%%%%%%%%%%%%%%%%%%%%%%%%%%%%%%%%%%%%%%%%%%%%%%%%%%%%%%%

Let us recall Oseledets' theorem in the two-dimensional area-preserving
case. A proof can be found in \cite{RefOseledets}.

\begin{theorem}
Let $f\in \Diff$. Then there exists a measurable
function $x\mapsto \lambda ^{+}(x)\geq 0$ such that 
\[
\lambda ^{+}(f,x)=\lim_{n\rightarrow +\infty }\frac{1}{n}\log \left\|
Df^{n}(x)\right\| 
\]
for a.e. $x\in M$. Moreover, if
$\mathcal{O}^{+}=\{x;\;\lambda^{+}(f,x)>0\}$ has positive measure then
for a.e. $x\in \mathcal{O}^{+}$ there is
a splitting $T_{x}M=E^{u}(x)\oplus E^{s}(x)$, depending measurably on $x$,
such that for $v\in T_{x}M-\{0\}$%
\begin{eqnarray*}
\lim_{n\rightarrow +\infty }\frac{1}{n}\log \left\| Df^{n}(x).v\right\|
&=&\left\{ 
\begin{tabular}{l}
$\lambda ^{+}(f,x)$ if $v\notin E^{s}(x)$ \\ 
$-\lambda ^{+}(f,x)$ if $v\in E^{s}(x)$%
\end{tabular}
\right. , \\
\lim_{n\rightarrow -\infty }\frac{1}{n}\log \left\| Df^{n}(x).v\right\|
&=&\left\{ 
\begin{tabular}{l}
$-\lambda ^{+}(f,x)$ if $v\notin E^{u}(x)$ \\ 
$\lambda ^{+}(f,x)$ if $v\in E^{u}(x)$%
\end{tabular}
\right. ,
\end{eqnarray*}
and
\[
\lim_{n\rightarrow \pm \infty }\frac{1}{n}\log \sin \measuredangle
(E^{u}(f^{n}x),E^{s}(f^{n}x))=0. 
\]
\end{theorem}

\begin{remark}
One can also define the lower Lyapunov exponent as 
\[
\lambda ^{-}(f,x)=\lim_{n\rightarrow -\infty }\frac{1}{n}\log \left\|
Df_{x}^{n}\right\| . 
\]
It satisfies $\lambda ^{-}(f,x)=-\lambda ^{+}(f,x)$.
\end{remark}

The \emph{integrated Lyapunov exponent} is defined by
$$ \mathrm{LE}(f)=\int_{M}\lambda ^{+}(f,x)d\mu(x) \, . $$

\begin{proposition}
\label{semicont}
The function $\mathrm{LE}:\Diff \to [0,\infty )$,
is upper semicontinuous. Besides, it is given by
$$
\mathrm{LE}(f)=\inf_{n \geq 1} \frac{1}{n}
\int_{M}\log \| Df^n \| d\mu \, .
$$
\end{proposition}

\begin{proof}
Given $f$, let $a_{n}(f)=\int_{M}\log \left\|
Df^{n}\right\| d\mu $. Then the sequence $(a_{n}(f))$ is subadditive, that is, $%
a_{n+m}\leq a_{n}+a_{m}$ for every $m,n$. Therefore 
\[
\mathrm{LE}(f)=\lim_{n\rightarrow +\infty }\frac{a_{n}(f)}{n}=\inf_{n\geq 1}%
\frac{a_{n}(f)}{n}\, .
\]
Thus LE$(\cdot)$ is the infimum of a sequence of continuous functions
and therefore it is upper semicontinuous.
\end{proof}

%%%%%%%%%%%%%%%%%%%%%%%%%%%%%%%%%%%%%%%%%%%%%%%%%%%%%%%%%%%%%%
\subsection{Avoiding periodic points and hyperbolic sets}
%%%%%%%%%%%%%%%%%%%%%%%%%%%%%%%%%%%%%%%%%%%%%%%%%%%%%%%%%%%%%%

We call $\mathcal{R}\subset \Diff$ a residual
subset if it contains a countable intersection of open dense sets. The set $%
\Diff$ is a Baire space, that is, every residual
subset is dense. A property is said to be \emph{generic} if it holds
in a residual set.

We say that a measure preserving transformation is \emph{aperiodic}
if the set of its periodic points has zero measure.

\begin{proposition}
For a dense set of $f\in \Diff$, the following holds:
Every hyperbolic set for $f$ has measure $0$ or $1$
and $f$ is aperiodic.
\end{proposition}

\begin{proof}
Let $f_0 \in \Diff$; we will find $f$ close to $f_0$ with the required properties.
Take a $C^2$ diffeomorphism $f_1 \in \Diff$ $C^1$-close to $f_0$
(it exists, as proved by Zehnder \cite{Zehnder}).
Using Robinson's conservative version of
Kupka-Smale theorem \cite[theorem 1.B.i]{Robinson}, we 
find a $C^2$ diffeomorphism $f$, $C^2$-close to $f_1$,
with countably many periodic points.
If $f$ is Anosov then we are done.
Otherwise,  we use that the hyperbolic sets for a $C^2$ non-Anosov diffeomorphism have
zero measure (this is a folklore theorem; for a
proof for basic sets see~\cite{Bowen}).
\end{proof}

It follows from proposition \ref{semicont} that the set
$\left\{ f\in \Diff;\,\mathrm{LE}(f)=0\right\} $
is a countable intersection of open subsets 
$\left\{ f;\,\mathrm{LE}(f)<\frac{1}{k} \right\} $ of $\Diff$.
Hence to prove theorem A we only need to show the following:

\begin{proposition}
\label{teorema}
Let $f\in \Diff$ be aperiodic and such that every
hyperbolic set for $f$ has zero measure. Let $\UU$
be a neighborhood of $f$ in $\Diff$
and let $\delta >0$. Then there exists $g\in \UU$
such that $\mathrm{LE}(g)<\delta $.
\end{proposition}

Proposition \ref{teorema} also has the following consequence:
\begin{proposition}
The functions 
$\mathrm{LE}$ and $h_{\mu}:\Diff\rightarrow \lbrack 0,\infty )$
are not continuous.
\end{proposition}

\begin{proof}
Any compact surface supports non-Anosov area preserving diffeomorphisms
with positive metric entropy, see Katok \cite{Katok}.
These diffeomorphisms, by proposition~\ref{teorema} and
Ruelle's inequality $h_\mu \leq \mathrm{LE}$, are points of discontinuity
of both functions.
\end{proof}

%%%%%%%%%%%%%%%%%%%%%%%%%%%%%%%%%%%%%%%%%%%%%%%%%%%%%%%%%%%
\subsection{Fixing coordinates, metrics and neighborhoods}
%%%%%%%%%%%%%%%%%%%%%%%%%%%%%%%%%%%%%%%%%%%%%%%%%%%%%%%%%%%

Now we establish some notation to be used until the end of section~\ref{s.prova teor A}.

Darboux's theorem (see \cite{Arnold}, for instance) gives that for each 
$x\in M$ there is an open set $V\ni x$ and a $C^{\infty }$ diffeomorphism 
$\varphi :V\rightarrow \varphi (V)\subset \RR ^{2}$ such that the induced
measure $\varphi _{\ast }(\mu )$ coincides with the usual Lebesgue measure
in $\varphi (V)\subset \RR ^{2}$. Taking a finite cover of $M$ by such
domains, we obtain an atlas 
$$
\mathcal{A}^{\ast }=\{\varphi _{i}:V_{i}^{\ast}\rightarrow \RR ^{2}, 
\; i=1,2, \ldots, \alpha \}.
$$
For technical reasons let us take open sets 
$V_{i}\subset M$ with $\overline{V_{i}}\subset V_{i}^{\ast }$ such that 
$\bigcup V_{i}=M$. Restricting the charts $\varphi _{i}$ to these smaller
domains we obtain another atlas $\mathcal{A}=\{\varphi _{i}:V_{i}\rightarrow 
\RR ^{2}\}$. 
We will also suppose that $\mu(\partial V_i) = 0 = \mu(\partial V_i^\ast)$
for each~$i$.

For each point $x \in M$, let 
$i(x) = \min \{i; \; V_i \ni x \}$.
We define a norm $\| \cdot \| = \| \cdot \|_x $ on each tangent space $T_x M$
by $\|v\| = \| D\varphi_{i(x)}(v) \|$.
The Riemannian metric on $M$ will not be used.

If $A: T_x M \to T_y M$ is a linear map, the norm
$\|A\|$ is then defined in the usual way:
$$
\|A\| = \sup_{0 \neq v \in T_x M} \frac{\| Av \|}{\| v \|}.
$$
Using the charts $\varphi_{i(x)}$ and $\varphi_{i(y)}$,
we may view $A$ as a linear map $\RR^2 \to \RR^2$.
This permits us, for example, to speak of the distance
$\|A  - B \|$ between two linear maps 
$A : T_{x_1}M \to T_{x_2}M$ and $B : T_{x_3}M \to T_{x_4}M$
whose base points are different.
Precisely, we define
$$
\|A - B \| = \| D_2 A D_1^{-1} - D_4 B D_3^{-1} \|, \quad
\text{where $D_j = (D\varphi_{i(x_j)})_{x_j}$.}
$$

If $x\in M$ and $r>0$ is small, we define
$$
B(x,r) = \varphi_{i(x)}^{-1} \big( B(\varphi_{i(x)}(x),r) \big).
$$
We will always assume that $r$ is small enough so that
$\overline{B(x,r)} \subset V_i^{\ast}$.
The sets $B(x,r)$ will be called \emph{disks}.

Given $\eps_0 > 0$, we define the
\emph{$\eps_0$-basic neighborhood $\UU (\mathrm{id}, \eps_0)$ of the identity} as the
set of diffeomorphisms $h \in \Diff$ such that we have
$h(V_i) \subset V_i^{\ast }$ for every $i$,
$h(x) \in B(x,\eps_0)$ and $\| Dh_{x} - I  \| < \eps_0$
for every $x \in M$.
Given $f \in \Diff$, we define the
\emph{$\eps_0$-basic neighborhood $\UU (f, \eps_0)$ of $f$} as the
the set of $g \in \Diff$ such that
$g \circ f^{-1} \in \UU(\mathrm{id}, \eps_0)$ or
$f^{-1} \circ g \in \UU(\mathrm{id}, \eps_0)$.
Notice that every neighborhood of $f$ in $\Diff$ contains a basic neighborhood.

%%%%%%%%%%%%%%%%%%%%%%%%%%%%%%%%%%%%%%%%%%%%%%%%%%%%%%%%%%%%%%%%%%
%%%%%%%%%%%%%%%%%%%%%%%%%%%%%%%%%%%%%%%%%%%%%%%%%%%%%%%%%%%%%%%%%%
\section{Construction of perturbations along an orbit segment}
%%%%%%%%%%%%%%%%%%%%%%%%%%%%%%%%%%%%%%%%%%%%%%%%%%%%%%%%%%%%%%%%%%
%%%%%%%%%%%%%%%%%%%%%%%%%%%%%%%%%%%%%%%%%%%%%%%%%%%%%%%%%%%%%%%%%%

The aim of this section is to prove the lemma below. In the next section we
will deduce theorem A from it.

\begin{main lemma}
Let $f\in \Diff$ be aperiodic and such that every 
hyperbolic set has zero measure, $\UU$ be a neighborhood of $f$ in
$\Diff$, $\delta >0$ and $0<k<1$. Then there exists a measurable integer
function $N:M\rightarrow \NN$ with the following properties:
For a.e. $x\in M$ and every integer $n\geq N(x)$ there exists $r=r(x,n)$ 
such that for every disk $U=B_{r^{\prime}}(x)$ with $0<r^{\prime}<r$,
there exist $g\in \UU $ and a compact set $K\subset U$ such that:
\begin{itemize}
\item[(i)] $g$ equals $f$ outside the set
$\bigsqcup_{j=0}^{n-1}f^{j}(\overline{U})$ and the iterates 
$f^{j}(\overline{U}),\, 0 \leq j \leq n-1$ are two-by-two
disjoint,

\item[(ii)] $\frac{\mu (K)}{\mu (U)}>k$,

\item[(iii)] if $y\in K$ then $\frac{1}{n}\log \left\|
Dg_{y}^{n}\right\| <\delta $.

\end{itemize}
\end{main lemma}

We will now outline some ideas in the proof of this lemma.
We perturb the derivatives $Df_{f^j x}$ of $f$ along the orbit segment.
These perturbations are constructed in such a way that their product
has small (i.e., not exponentially large) norm. To hinder the growth
of these products, we send the expanding Oseledets direction in the
contracting one. This is possible in the absence of uniform hyperbolicity.
Once the linear perturbations are constructed, we must find an area preserving
diffeomorphism $g$ close to $f$ having approximately
the assigned derivatives. To guarantee that the diffeomorphism
$g$ exists and has the stated properties, some care is needed
in the choice of the perturbations of $Df_{f^j x}$.

%%%%%%%%%%%%%%%%%%%%%%%%%%%%%%%%%%%%%%%%%%%%%%%
\subsection{Realizable sequences}
%%%%%%%%%%%%%%%%%%%%%%%%%%%%%%%%%%%%%%%%%%%%%%%

As we mentioned, to prove the Main Lemma we will construct perturbations
$L_{j}$ of the linear maps $Df_{f^j x}$. These perturbations will
be required to have the following property:

\begin{definition}
Given $f\in \Diff$, a basic neighborhood
$\UU =\UU (f,\varepsilon _{0})$, $0<k<1$ and a
non-periodic point $x \in M$, a sequence of (area form preserving) linear maps 
\[
T_{x}M\stackrel{L_{0}}{\longrightarrow }T_{fx}M\stackrel{L_{1}}{%
\longrightarrow }\cdots \stackrel{L_{n-1}}{\longrightarrow }T_{f^{n}x}M 
\]
is called a $(k,\UU )$\emph{-realizable sequence} of length $n$ at $x$
if the following holds: For every $\gamma >0$ there is $r>0$ such that if 
$U\subset B_{r}(x)$ is a non-empty open set then there are $g\in \UU $
and a compact set $K\subset U$ such that
\begin{itemize}
\item[(i)] $g$ equals $f$ outside the set
$\bigsqcup_{j=0}^{n-1}f^{j}(\overline{U})$ and the iterates 
$f^{j}(\overline{B_{r}(x)}),\, 0 \leq j \leq n-1$ are two-by-two
disjoint,

\item[(ii)] $\frac{\mu (K)}{\mu (U)}>k,$

\item[(iii)] if $y\in K$\ then $\left\| Dg_{g^{j}y}-L_{j}\right\| <\gamma 
$ for every $j$.

\end{itemize}
\end{definition}

In the following lemma we give some useful properties of realizable sequences:

\begin{lemma}\label{R}
\begin{enumerate}
\item The sequence $\{Df_{x},\ldots ,Df_{f^{n-1}(x)}\}$ is 
$(k,\UU )$-realizable for every $k$ and $\UU $.

\item Let $k$, $k_1$, $k_2 \in (0,1)$ be such that 
$1-k=(1-k_1)+(1-k_2)$.
If $\{L_{0},\ldots ,L_{n-1}\}$ is a $(k_1,\UU )$-realizable sequence at $x$
and $\{L_{n},\ldots,L_{n+m-1}\}$ is a $(k_2,\UU )$-realizable sequence at $f^{n}(x)$
then $\{L_{0},\ldots ,L_{n+m-1}\}$ is $(k,\UU )$-realizable at $x$.

\item Let $\{L_j : T_{f^j x}M \to T_{f^{j+1}x} M\}$ be a sequence of linear maps
at a non-periodic point $x$.
To prove that the sequence $\{L_j\}$ is $(k,\UU )$-realizable
we only need to check the conditions for open sets $U$ that are
disks, $U=B_{r_{0}}(y)\subset B_{r}(x) \subset V$.
\end{enumerate}
\end{lemma}

\begin{proof}
For the first property, just take $g=f$.

For property (2), take $\gamma>0$. Let $r_1>0$ be the radius associated to the first
($(k_1,\UU )$-realizable) sequence, and $r_2$ the radius associated
to the second sequence.
Let $0<r<r_1$ be such that $f^n(B_r(x)) \subset B(f^n(x),r_2)$. Given an open set
$U \subset B_r(x)$, the realizability of the first sequence gives us a diffeomorphism
$g_1 \in \UU $ and a set $K_1 \subset U$. Analogously, from the open set
$f^n(U) \subset B(f^n(x),r_2)$ we can find $g_2 \in \UU $ and $K_2 \subset f^n(U)$.
Then define a diffeomorphism $g$ as 
$g=g_1$ inside $U,\ldots,f^{(n-1)}(U)$,
$g=g_2$ inside $f^n(U),\ldots,f^{(n+m-1)}(U)$, $g=f$ elsewhere.
Define also a compact set $K = K_1 \cap g^{-n}(K_2) \subset U$.
Then one can check that $g$ and $K$ satisfy the required properties.

Now let us prove (3). 
Let $\{L_j\}$ be a sequence at the point $x$ and suppose that the conditions of
realizability are satisfied for open sets $U$ that are disks.
That is, given $\gamma>0$, there is $r>0$ such that for every
disk $U \subset B(x,r)$ there are $g$ and $K$ verifying conditions (i)-(iii)
of the definition of a realizable sequence.
Fix the chart $(\varphi ,V) 
=(\varphi_{i(x)}, V_{i(x)}) \in \mathcal{A}$
and take any open set $U\subset B_r(x)$.
By Vitali's covering lemma (see, for instance, \cite{McShane}),
there is a countable family of disjoint closed disks in $\RR ^2$
covering $\varphi_i^{-1}(U)$ mod 0.
Thus we can find a finite family of disks $U_i=B_{r_{i}}(y_{i})\subset U$
with disjoint closures such that $\mu \left(
U-\bigsqcup_{i}B_{r_{i}}(y_{i})\right) $ is as small as we please.
For each disk $U_i$ there are, by hypothesis,
a perturbation $g_{i}\in \UU (f,\varepsilon _{0})$ and a set $K_i \subset U_i$
with the properties (i)-(iii) of relizability.
Let $K=\bigcup K_i$ and define $g$ as equal $g_i$ in each $f^j(U_i)$.
Since the latter sets are disjoint, $g$ is well-defined.
Moreover, $g\in \UU (f,\varepsilon _{0})$ and the pair $(g,K)$ satisfies
the required properties (i)-(iii).
\end{proof}

We will denote by $R_{\theta}$ the rotation of angle $\theta $,
\[
R_{\theta }=\left( 
\begin{array}{cc}
\cos \theta & -\sin \theta \\ 
\sin \theta & \cos \theta
\end{array}
\right) . 
\]
The simple lemma below is the basic tool that will be used to construct all
our realizable sequences.

\begin{lemma}\label{basico}
Let $\varepsilon _{1}>0$ and $0<k<1$. Then there exists 
$\alpha _{0}>0$ with the following properties: If $|\alpha |\leq
\alpha _{0}$ and $r>0$, then there exists a $C^{1}$ area preserving
diffeomorphism $h:\RR ^{2}\rightarrow \RR ^{2}$ such that
\begin{itemize}

\item[(i)] $\left| z\right| \geq r\Rightarrow h(z)=z$,

\item[(ii)] $\left| z\right| \leq \sqrt{k} r\Rightarrow h(z)=R_{\alpha }(z)$,

\item[(iii)] $|h(z)|=|z|\quad \forall z$,

\item[(iv)] $\left| h(z)-z\right| \leq \alpha r\quad \forall z$,

\item[(v)] $\left\| Dh_{z}-I\right\| \leq \varepsilon _{1}\quad \forall z$.
\end{itemize}
\end{lemma}

\begin{proof}
First suppose $r=1$. Let $F:\RR \rightarrow \RR $
be a $C^{\infty }$ function such that
$F(t)=1$ for $t\leq \sqrt{k}$, $F(t)=0$ for $t\geq 1$ and
$0\leq -F^{\prime }(t)\leq\frac{2}{1-\sqrt{k} }$ .
Define $h$ by $h(z)=R_{\alpha F(\left| z\right| )}(z)$. Then $h$ is an
area preserving diffeomorphism satisfying properties (i)-(iv). If $\alpha
_{0}$ is small, (v) will hold for $\left| z\right| \leq \sqrt{k} .$ We still have
to check (v) in the annulus $\sqrt{k} \leq \left| z\right| \leq 1,$ where we can
take polar coordinates $(\rho ,\theta )$. The mapping $h$ takes the form $%
(\rho ,\theta )\mapsto (\rho ,\theta +\alpha F(\rho ))$.
The Jacobian matrix of $h$ is
$\left( 
\begin{array}{cc}
1 & 0 \\ 
\alpha F^{\prime }(\rho ) & 1
\end{array}
\right) $ and is close to identity if $\alpha _{0}$ is small enough.

Now we claim that the same $\alpha _{0}$ will work for any $r>0$.
Indeed, if $h_{1}$ is the diffeomorphism associated to $r=1$
constructed above, let $h(z)=rh_{1}(r^{-1}z)$.
Then (i)-(iv) obviously hold for $h$ and, since 
$Dh_{z}=D(h_{1})_{r^{-1}z}$, (v) also holds.
\end{proof}

The two lemmas below will be used to construct realizable sequences.
Given $x \in M$ and $\theta \in \RR$, consider the chart
$\varphi =\varphi_{i(x)} : V_{i(x)} \to \RR ^2$ and the linear map
$(D\varphi_x)^{-1} R_\theta D\varphi_x :T_{x}M \hookleftarrow$.
We shall call this map \emph{the rotation of angle $\theta$ at $x$}
and denote it also by $R_\theta$.

\begin{lemma} \label{tipo1}
Given $f\in \Diff$,
$\UU =\UU (f,\varepsilon _0)$ and $0<k<1$,
there is $\alpha_{1} >0$ with the following properties:
Suppose that $x\in M$ is not periodic and $| \theta | \leq \alpha_{1}$.
Then $\{Df_{x} R_{\theta}\}$ and $\{R_{\theta} Df_{x} \}$ are
$(k,\UU )$-realizable sequences of length $1$ at $x$.
\end{lemma}

\begin{proof}
We will prove that if $|\theta|$ is small then the sequence
$\{Df_{x} R_{\theta}\}$ is realizable; for the other sequence the proof is similar.
Let $r>0$ be small.
By lemma \ref{R}.3, we need only to construct perturbations supported in disks 
$U=B_{r^{\prime}}(y) \subset B_r(x)$.
Now apply lemma \ref{basico} to find, for each small angle $\theta$, a diffeomorphism
$g$ and a disk $K \subset U$ with the required properties.
Since $r$, and thus $r^{\prime}$, is small, $g$ is near $f$.
\end{proof}

\begin{lemma} \label{tipo2}
Given $f\in \Diff$,
$\UU =\UU (f,\varepsilon _{0})$ and $0<k<1$,
there is $\alpha_{2} >0$ with the following properties:
Suppose that $x\in M$ is not periodic, $m \geq 2$,
$R_{\theta_{0}}:T_{x}M\hookleftarrow $ and
$R_{\theta_{1}}:T_{f^{m}(x)}M\hookleftarrow $ are rotations such that
$|\theta_{i}| \leq \alpha_{2}$. Then 
$$
\{Df_{x} R_{\theta_{0}},Df_{fx}, \dots ,
Df_{f^{m-2}x}, R_{\theta_{1}} Df_{f^{m-1}x}\}
$$
is a $(k,\UU)$-realizable sequence of length $m$ at $x$.
\end{lemma}

\begin{proof}
Let $k_0$ be such that $1-k_0<\frac{1}{2}(1-k)$. Apply lemma \ref{tipo1} with 
$k_0$ in the place of $k$ to obtain $\alpha_2$. 
Now, if $|\theta_{0}|$, $|\theta_{1}| \leq \alpha_{2}$ then 
$\{Df_{x}R_{\theta_{0}}\}$ and $\{R_{\theta_{1}}Df_{f^{m-1}x}\}$ are 
$(k_0,\UU )$-realizable sequences of length $1$. By items 1 and 2 of lemma \ref{R},
$\{Df_{x}R_{\theta_{0}},Df_{fx}, \dots , Df_{f^{m-2}x},R_{\theta_{1}}Df_{f^{m}x}\}$ is
a $(k,\UU)$-realizable sequence.
\end{proof}

The lemma above is somewhat weak; we cannot use
an arbitrary (say, $m$) number of rotations instead of just two.
This difficulty will be overcome in the next section.

%%%%%%%%%%%%%%%%%%%%%%%%%%%%%%%%%%%%%%%%%%%%%%%
\subsection{Nested rotations}
%%%%%%%%%%%%%%%%%%%%%%%%%%%%%%%%%%%%%%%%%%%%%%%

Now we will deal with linear area preserving transformations
with an invariant ellipse. As it will be shown in lemma \ref{tipo3},
suitable sequences of such nested rotations are realizable.
The key point there is that those sequences can be
arbitrarily long while $k$ is kept controlled.

Our ellipses will always be filled ones. An ellipse
$\BB \subset \RR ^{2}$ has 
\emph{eccentricity} $E$ if it is the image of a disk under a
transformation $L\in \sldoisr$ with $\| L\| =E$. (This is not
quite the usual definition). Thus:
$$
E=\sqrt{\frac{\mbox{major axis}}{\mbox{minor axis}}} \, .
$$
 
The following lemma is a slight generalization of lemma \ref{basico}.

\begin{lemma}\label{rotaelip}
Let $\varepsilon _{1}>0$, $0<k<1$ and $E\geq 1$.
Then there exists $\varepsilon >0$ with the following properties:
Let $\BB \subset \RR ^{2}$ be an ellipse centered 
at the origin, of eccentricity $\leq E$
and diameter $\leq \varepsilon _{1}$, and 
let $L\in \sldoisr$ be a linear map 
with $\left\| L-I\right\| <\varepsilon $ and
preserving $\BB $. Then there exists a $C^{1}$ area preserving
diffeomorphism $h:\RR ^{2}\rightarrow \RR ^{2}$ such that if
$\BB ^{\prime }$ is the (smaller) ellipse 
$\BB ^{\prime }=\sqrt{k}\BB $ then
\begin{itemize}
\item[(i)] $z\notin \BB \Rightarrow h(z)=z$,

\item[(ii)] $z\in \BB ^{\prime }\Rightarrow h(z)=L(z)$,

\item[(iii)] $h$ preserves all ellipses of the form $t\BB $, $t>0$,

\item[(iv)] $\left| h(z)-z\right| \leq \varepsilon _{1}\quad \forall z$,

\item[(v)] $\left\| Dh_{z}-I\right\| \leq \varepsilon _{1}\quad \forall z$.

\end{itemize}
\end{lemma}

\begin{proof}
There is $M \in \sldoisr$ with $\|M\| = \|M^{-1}\| = E$
such that $M(\BB )$ is a disk.
We apply lemma \ref{basico} with $E^{-2}\varepsilon _{1}$
in the place of $\varepsilon _{1}$ to get $\alpha _{0}>0$.
If $L\in \sldoisr$ preserves $\BB $ and
$\left\| L-I\right\| <\varepsilon $ then
$\| MLM^{-1}-I\| <E^{2}\varepsilon $. Thus if $\varepsilon >0$
is chosen small enough then $MLM^{-1}$ is a rotation of angle $\alpha $
with $|\alpha |<\alpha _{0}$. Thus lemma \ref{basico} gives a
diffeomorphism $h_{0}$ and we define $h=M^{-1}h_{0}M$. The properties
(i)-(iii) and (v) are easily seen to hold for $h$; (iv) is assured 
if we suppose that $\varepsilon  $, and thus $|\alpha |$, is very small.
\end{proof}

The next lemma says that the image of a small ellipse by a $C^{1}$
diffeomorphism is approximately an ellipse.

\begin{lemma}\label{elao}
Let $h:\RR ^{2}\rightarrow \RR ^{2}$ be a $C^{1}$ diffeomorphism with
$h(0)=0$. Let $M=Dh_{0}:\RR ^{2}\rightarrow \RR ^{2}$. Then, given
$\eta >0$ and $E>1$ there exists $r>0$ such that if $\BB (y)\subset
B_{r}(0)$ is an ellipse centered in $y$ with eccentricity $\leq E$ and 
$\BB (0)=\BB (y)-y$ is the translated ellipse centered in $0$,
then we have 
\[
(1-\eta )M\BB (0)+h(y)\subset h(\BB (y))\subset (1+\eta )M\BB (0)+h(y). 
\]
\end{lemma}

\begin{proof}
Let $g:\RR ^{2}\rightarrow \RR ^{2}$ be such
that $h=Mg$. Since $g$ is $C^{1}$ smooth and $Dg_{0}=I$, we have 
\begin{equation}
g(z)-g(y)=z-y+\xi (z,y)  \tag{$\ast $}
\end{equation}
where 
\[
\lim_{(z,y)\rightarrow (0,0)}\frac{\xi (z,y)}{\left| z-y\right| }=0. 
\]
($| \cdot |$ indicates the euclidean norm in $\RR ^{2}$). Choose 
$r>0$ such that
$| z| ,| y| <r\Rightarrow | \xi
(z,y)| \leq E^{-1}\eta | z-y| $. Now let
$\BB (y)\subset B_{r}(0)$ be an ellipse with axes $2a$ and $2b$;
$1\leq \frac{b}{a}\leq E$.
If $z\in \partial \BB (y)$ then
$|z-y| \leq b$ and
$| \xi (z,y)| \leq E^{-1}\eta b\leq \eta a$. Thus equation ($\ast $)
gives that
\[
g(\partial \BB (y))-g(y)\subset V_{\eta a}(\partial \BB (0)), 
\]
where $V_{\varepsilon }(\cdot )$ denotes $\varepsilon $-neighborhood. To avoid confusion,
we will temporarily denote difference of sets by the symbol $\setminus $.
We have the geometrical property 
\[
V_{\eta a}(\partial \BB (0))\subset
(1+\eta )\BB (0) \setminus  (1-\eta )\BB (0).
\]
Therefore 
\[
g(\partial \BB (y))-g(y)\subset
(1+\eta )\BB (0) \setminus  (1-\eta )\BB (0).
\]
Applying the linear mapping $M$, we have 
\[
h(\partial \BB (y))-h(y)\subset
(1+\eta )M\BB (0) \setminus  (1-\eta )M\BB (0).
\]
The lemma follows now from standard topological arguments.
\end{proof}

Now we use the above material to give another construction
of realizable sequences.
In the next lemma,
$\mathbb{D} \subset T_{x}M$ will denote
the unit disk $\{v; \; \|v\|<1\}$.

\begin{lemma}\label{tipo3}
Given $f\in \Diff$, $\UU =\UU(f,\varepsilon _{0})$, $0<k<1$ and $E>1$,
there is $\varepsilon >0$ with the
following properties: Suppose that $x\in M$ is not periodic and there is
$n\in \NN$ such that $\left\| Df_{x}^{j}\right\| \leq E$ for 
$j=1,\ldots,n$. If 
\[
T_{x}M\stackrel{L_{0}}{\longrightarrow }T_{fx}M\stackrel{L_{1}}{%
\longrightarrow }\cdots \stackrel{L_{n-1}}{\longrightarrow }T_{f^{n}x}M 
\]
are linear maps such that for every $j=1,\ldots ,n$ we have:
\begin{itemize}
\item[(i)] $L_{j-1}\ldots L_{0}(\mathbb{D})=Df_{x}^{j}(\mathbb{D})$ and

\item[(ii)] $\left\| L_{j}-Df_{f^{j}x}\right\| <\varepsilon$,

\end{itemize}
then $\{L_{0},\ldots ,L_{n-1}\}$ is $(k,\UU)$-realizable at $x$.
\end{lemma}

\begin{proof}
Let $f$, $\eps_0$, $k$ and $E$ be given.
Let $\eps>0$ be given by lemma~\ref{rotaelip}, depending on
$\eps_1 = \eps_0$, $k$ and $E$.

Now let $x$, $n$ and $\{L_0,\ldots,L_{n-1}\}$ be as in the statement.
We must prove that the sequence $\{L_0,\ldots,L_{n-1}\}$ is $(k,\UU)$-realizable;
so let $\gamma > 0$ be given.

We will consider the charts $\varphi_{i(f^j x)} : V_{i(f^j x)} \to \RR^2$.
To simplify notation, we write $\varphi_j : V_j \to \RR^2$ instead.

Let $r_0>0$ be such that, for each $j=0,1\ldots,n$, we have:
\begin{itemize}
\item[(1)] $f^{j}( \overline{B_{r_0}(x)})\subset V_{i(f^j(x))}^{\ast },$

\item[(2)] the sets $f^{j}(\overline{B_{r_0}(x)})$ are two-by-two
disjoint,

\item[(3)] for every $z\in f^{j}(B_{r_0}(x))$ we have $\left\|
Df_{z}-Df_{f^{j}x}\right\| <\gamma E^{-2}$,

\item[(4)] $\|Df_z^j\| < 2E$ for every $z \in B_{r_0}(x)$.

\end{itemize}

Using the charts, we can translate the problem to $\RR ^{2}$.
Let $f_{j}$ be the expression of $f$ in charts in the neighborhood of $f^{j}(x)$,
that is, $f_{j}=\varphi_{j+1}\circ f\circ \varphi _{j}^{-1}$.
To simplify notations, we suppose that $\varphi _{j}(f^{j}(x))=0$.

Let $\tilde{L}_j$ be the expression of $L_j$ in charts, that is,
$$
\tilde{L}_j = (D\varphi_{j+1})_{f^{j+1} x} \cdot L_j \cdot [(D\varphi_j)_{f^j x}]^{-1}.
$$
Let $M_{j}=(Df_{j})_{0}\in \sldoisr$. By hypothesis, we have
$\left\| M_{j-1}\ldots M_{0}\right\| \leq E$ for each $j$.
Let $R_{j} \in \sldoisr$ be such that $\tilde{L}_{j} = M_{j} R_{j}$. 
Since $R_j$ preserves an ellipse of eccentricity $\leq E$, we have
$\|R_j\| \leq E^2$.

Take $k_{0}$ such that $k<k_{0}<1$ and let $0<\tau <1$ be such that
$\tau ^{4n}k_{0}>k$. Take $\eta>0$ such that
$\tau < 1-\eta < 1+\eta < \tau^{-1}$.
Using lemma~\ref{elao} $n$ times, we find $0<r_1<r_0$ with the following properties:
If $\BB \subset B_{r_1}(0)$ is an ellipse centered at a point $y$ and 
with eccentricity $\leq E$ then
$$
(1-\eta) M_j (\BB - y ) \subset f_j(\BB) - f_j(y) \subset (1+\eta) M_j (\BB - y )
\quad \text{for each $j=0,\ldots,n-1$.}
$$
Finally, let $r = r_1 /3E$. 
%such that $f^j(B_r(x)) \subset B_{r_1}(f^j(x))$
%(which is equivalent to $f_{j-1} \cdots f_0 (B_r(0)) \subset B_{r_1}(0)$)
%for each $j=0,\ldots,n-1$.
This is the $r$ of the definition of realizable sequences.

By lemma \ref{R}.3, in order to prove realizability we may restrict
ourselves to sets $U$ that are disks,
so let
$U = B_{r'}(x') \subset B_r(x)$, for some $x'$ and $r'$.

Let, for $0< t \leq 1$, 
$\BB _{t}^{0}=B_{tr'}(x')$.
Let $\BB _{t}^{0}-x'$ denote the translate of $\BB _{t}^{0}$
centered at the origin. Define, for $t>0$ and $j=1,\ldots ,n$, 
$$
\BB_t^j = M_{j-1}\cdots M_{0}(\BB _{t}^{0}-x')+
f_{j-1}\cdots f_{0}(x')\, .
$$
Then $\BB_t^j$ is an ellipse centered at
$f_{j-1}\cdots f_{0}(x')$ with eccentricity $\leq E$. 
We have $M_{j-1}\cdots M_{0} ( B_r(0) ) \subset B_{Er}(0)$ and,
by (4), $f_{j-1}\cdots f_{0} ( B_r(0) ) \subset B_{2Er}(0)$.
Therefore
$$
\BB_t^j \subset \BB_1^j \subset M_{j-1}\cdots M_{0} ( B_r(0) ) + f_{j-1}\cdots f_{0}(x') 
                        \subset B_{3Er}(0) \subset B_{r_1}(0).
$$
This permits us to apply lemma~\ref{elao} to those ellipses and get the property
$$
0 < t \leq 1, \  0 \leq j \leq n-1 \Longrightarrow
\BB_{\tau t}^{j+1} \subset f_j(\BB_t^j) \subset \BB _{\tau ^{-1}t}^{j+1} \, .
$$

By hypothesis (i), the linear map $R_{j}$ preserves the ellipses $\BB _{t}^{j}$.
So we apply lemma \ref{rotaelip} with $k_{0}$,
$\BB _{\tau ^{n}}^{j}$ and $R_{j}$ in the place of
$k$, $\BB $ and $L$. 
The lemma gives us a $C^{1}$ area preserving diffeomorphism 
$h_{j}:\RR ^{2}\rightarrow \RR ^{2}$ such that:
\begin{itemize}

\item[(5)]  $z\notin \BB _{\tau ^{n}}^{j}\Rightarrow h_{j}(z)=z$,

\item[(6)]  $z\in \overline{\BB}_{\tau ^{n}\sqrt{k_{0}}}^{j}
\Rightarrow h_{j}(z)=R_{j}(z)$,

\item[(7)]  $t \geq 0 \Rightarrow h_{j}(\BB _{t}^{j})=\BB _{t}^{j}$,

\item[(8)]  $| h_{j}(z)-z| <\varepsilon _{0}\quad \forall z$,

\item[(9)]  $\| D(h_{j})_{z}-I\| \leq \varepsilon_{0}\quad \forall z$.

\end{itemize}

For $0\leq j \leq n-1$, let $S_j = \varphi_j^{-1} (\{z;\; h_j(z)\neq z\})$.
Then we have
\begin{multline*}
S_j \subset \varphi_j^{-1} (\BB_{\tau^n}^j)
    \subset \varphi_j^{-1} \big[ f_{j-1} (\BB_{\tau^{n-1}}^{j-1}) \big]
    \subset \cdots \\
    \subset \varphi_j^{-1} \big[ f_{j-1}  \cdots  f_0 (\BB_{\tau^{n-j}}^0) \big] 
    =       f^j ( B_{\tau^{n-j} r'} (x')) 
    \subset f^j (U).
\end{multline*}
So, by (1) and $U \subset B_{r_0}(x)$, the sets $\overline{S_j}$ are disjoint.
This permits us to define a diffeomorphism $g\in \Diff$ 
as equal to $f$ in $S_0 \sqcup\cdots \sqcup S_{n-1}$, and 
equal to $\varphi_{j+1}^{-1}\circ g_{j}\circ \varphi_{j}$ in $S_j$,
where $g_{j}=f_{j}\circ h_{j}$.

Let us verify that $g$ is the desired perturbation. 
First of all, by (8) and (9), we have 
$f^{-1} \circ g \in \UU(\mathrm{id},\eps_0)$, that is,
$g \in \UU(f, \eps_0)$.

We abbreviate $f^{(0)} = \mathrm{id}$, $f^{(j)}=f_{j-1}\cdots f_{0}$ and
analogously for $g^{(j)}$.
Condition (i) in the definition of realizable sequences is easy to check:
We have seen that $f^j(U) \supset S_j$;
this means that $g_{j}$ equals $f_{j}$ outside $f^{j}(U)$.
Now define 
$$
K = \varphi_0^{-1} \big( \overline{\BB}^0_{\tau ^{2n}\sqrt{k_{0}}} \big)
=\overline{B}_{\tau ^{2n}\sqrt{k_{0}}}(x') \subset U \, .
$$
Since $\varphi_0$ takes $\mu$ to the area in $\RR ^{2}$, we can calculate
$$
\mu (U-K)= \pi  r^{\prime 2} - \pi \tau ^{4n}k_{0}r^{\prime 2}
=\left( 1-\tau ^{4n}k_{0}\right) \mu (U)<(1-k)\mu (U)\;.
$$
which is condition (ii) in the definition of realizable sequences.

Since $g_{j}=f_{j}\circ h_{j}$ and $h_{j}$ preserves $\BB _{t}^{j}$
ellipses, we have (by induction in $j$)
$$
g^{(j)} ( \overline{\BB}^0_t ) \subset \overline{\BB} ^j _{ \tau ^{-j}\sqrt{k_0} t}
\subset \overline{\BB} ^j _{\tau ^{-n}\sqrt{k_0} t}
$$
for every $j$ and $0<t<1$.
Setting $t=\tau ^{2n}\sqrt{k_{0}}$, we get
$$
g^j(K) \subset \varphi_j^{-1} \big( \overline{\BB}^j _{\tau ^{n}\sqrt{k_0} } \big).
$$

To check condition (iii), take $y \in K$. 
Let $\tilde{y} = \varphi_0^{-1}(y)$.
Then
$g^{(j)}(\tilde{y}) \in \overline{\BB} ^j _{\tau ^{n}\sqrt{k_0} }$
and so, by (6),
$\left( Dh_j \right)_{g^{(j)} (\tilde{y})} = R_j$.
Therefore
$$
\| Dg_{g^j y } - L_j \| =
\| (Dg_{j})_{g^{(j)} \tilde{y} } - \tilde{L}_{j} \| =
\| (Df_{j})_{h_{j}(g^{(j)} \tilde{y})}R_{j} - (Df_{j})_{0} R_{j} \| . 
$$
Using (3) and $\| R_{j} \| \leq E^2$ we get
$$
\| Dg_{g^j y } - L_j \| \leq
\|(Df_{j})_{h_{j}(g^{(j)} \tilde{y})} - (Df_{j})_{0} \| \;
\| R_{j} \|<
\gamma \, ,
$$
proving the third condition and thus the lemma.
\end{proof}

%%%%%%%%%%%%%%%%%%%%%%%%%%%%%%%%%%%%%%%%%%%%%%%%%%%%%%%%%%%%%%%
\subsection{Sending $E^{u}$ to $E^{s}$} \label{s:Eu em Es}
%%%%%%%%%%%%%%%%%%%%%%%%%%%%%%%%%%%%%%%%%%%%%%%%%%%%%%%%%%%%%%%

Here we use lemmas \ref{tipo1}, \ref{tipo2} e \ref{tipo3}
to construct realizable sequences that send the expanding Oseledets direction
in the contracting one; this is the content of lemma~\ref{Eu em Es} below.

First we define some notation that will also be used in section \ref{s:prova ML}.
Given a diffeomorphism $f \in \Diff$, let $\mathcal{O}(f)\subset M$
be the full measure set given by Oseledets' theorem.
Define the following $f$-invariant sets:
\begin{align*}
\mathcal{O}^{+}(f) &= \left\{ x\in \mathcal{O}(f);\;\lambda ^{+}(f,x)>0\right\},  \\
\mathcal{O}^{0}(f) &= \left\{ x\in \mathcal{O}(f);\;\lambda ^{+}(f,x)=0\right\} .
\end{align*}
Now we define for $x\in \mathcal{O}^{+}(f)$ and $m\geq 1$,
\[
\Delta (x,m) = \frac{\left\| Df_{x}^{m}\mid _{E^{s}(x)}\right\| }{\left\|
Df_{x}^{m}\mid _{E^{u}(x)}\right\| } 
\]
and define the set:
\[
\Gamma_m(f) = \left\{ x\in \mathcal{O}^{+}(f);\; \Delta(x,m) \geq 1/2 \right\}.
\]

In informal words, the set $\Gamma_m(f)$, for large $m$,
is the place where the lack of uniform hyperbolicity
appears, and where the Oseledets directions can be mixed.
More precisely, we have:
 
\begin{lemma} \label{Eu em Es}
Let $f\in\Diff$, $\UU = \UU(f,\eps_0)$ and $0<k<1$.
Then there is a positive integer $m$ such that for every
point $x \in \Gamma_m(f)$, there exists a $(k,\UU)$-realizable sequence
$\{L_0, L_1, \ldots, L_{m-1}\}$ at $x$ of length $m$ such that
$$
L_{m-1} \cdots L_1 L_0 (E^u(x)) = E^s (f^m(x)).
$$
\end{lemma}

For the proof of lemma~\ref{Eu em Es}, we will need two
simple linear-algebraic lemmas:

\begin{lemma} \label{linear1}
Given $\alpha_2 > 0$, there is $c>1$ with the following properties:
Given a linear transformation $A:\RR^2 \to \RR^2$ and unit vectors $s$, $u \in \RR^2$
such that $\| A(s) \| / \| A(u) \| > c$, there exists $\xi \in \RR^2 - \{0\}$ such that
$$
\ang(\xi, u) \leq \alpha_2 \quad \text{and} \quad \ang(A(\xi), A(s)) \leq \alpha_2.
$$
\end{lemma}

\begin{proof}
Let $a = \arcsin\alpha_2$ and $c=a^{-2}$.
Let $\xi = u + a s$. Since $\|u\|=\|s\|=1$, we have $\sin \ang(\xi,u) \leq a$.
Let $\hat{u} = A(u)/\| A(u) \|$ and $\hat{s} = A(s)/\| A(s) \|$.
The vector $A(\xi) = A(u) + a A(s)$ is colinear to
$$
\Big( a^{-1} \frac{\|A(u)\|}{\| A(s) \|} \Big) \hat{u} + \hat{s}.
$$
Therefore 
$\sin \ang(A(\xi),\hat{s}) \leq a^{-1} \frac{\|A(u)\|}{\| A(s) \|} < a$.
\end{proof}

\begin{lemma} \label{linear2}
Given $\alpha_1 > 0$ and $\hat{c}>1$, there exists $E>1$ with the following properties:
Let $A\in \sldoisr$.
If there exists a pair of unit vectors $s$, $u \in \RR^2$
such that $\ang(u,s) \geq \alpha_1$, $\ang(A(u),A(s)) \geq \alpha_1$ and
$$
\frac{1}{\hat{c}} \leq \frac{\| A(s) \|}{\| A(u) \|} \leq \hat{c},
$$
then $\|A\| \leq E$.
\end{lemma}

\begin{proof}
By substituting $A$ by $B_2^{-1} A B_1$, where $B_1$ and $B_2$ are area-preserving changes 
of coordinates whose norms are bounded by some function of $\alpha_1$, we may suppose that
$s \perp u$ and $A(s) \perp A(u)$.
Now we have $\| A(s) \| \cdot \| A(u) \| = \det A = 1$, so both $\| A(s) \|$
and $\| A(u) \|$ are bounded by $\hat{c}^{1/2}$.
\end{proof}

\begin{proof}[Proof of lemma ~\ref{Eu em Es}]
We first define constants $k_0$, $C$, $\alpha_1$, $\alpha_2$, $E$, $\eps$, $C$, $\beta$ and $m$.

Let $k_0 \in (0, k)$.
Let $C = \sup \|Df^{\pm 1}\|$.
Let $\alpha_1>0$ and $\alpha_2>0$, depending on $f$, $\eps_0$ and $k_0$ (in the place of $k$),
be given by lemmas~\ref{tipo1} and \ref{tipo2}, respectively. 
Let $c$, depending on $\alpha_2$, be given by lemma~\ref{linear1}.
We assume that $c > C^2$.

Let $E>1$, depending on $\alpha_1$ and $\hat{c}=2c^2$, be given by lemma~\ref{linear2}.

Let $\eps> 0$, depending on $f$, $\eps_0$, $k$ and $E$, be given by lemma~\ref{tipo3}.

Choose $\beta>0$ such that if $|\theta | \leq \beta$ then the rotation $R_{\theta}$ is
close to the identity, $\| R_{\theta} - I \| < C^{-1}E^{-2}\varepsilon$.
Let $m$ be the least integer satisfying $m \geq \frac{2\pi}{\beta}$.

Fix $x\in \Gamma_m(f)$. The rest of proof is divided in three cases.

\subsubsection*{First case}
Suppose that the following condition holds:
\begin{equation}\tag{I}
\text{there exists $j_0 \in \{0,1,\ldots,m-1\}$ such that }
\ang(E^u(f^{j_0}(x)), E^s(f^{j_0}(x)) < \alpha_1.
\end{equation}

By lemma \ref{tipo1}, for every $y \in M$ and every
$|\theta| < \alpha_1$ and every rotation $R_\theta$ at $y$, the sequence 
$\{Df_y R_\theta\}$ of length $1$ at $y$ is $(k,\UU)$-realizable.
We use this fact with $y=f^{j_0}(x)$ and $\theta = \pm \ang(E^u(y), E^s(y))$
such that $R_\theta (E^u(y)) = E^s(y)$.
So $\{Df_{f^{j_0}(x)} R_\theta\}$ is a $(k_0,\UU)$-realizable sequence
of length $1$ at $y$.
By lemma \ref{basico}, items 1 and 2, the sequence
$$
\{L_0, \ldots, L_{m-1} \} =
\{Df_x, \ldots, Df_{f^{j_0-1}(x)}, Df_{f^{j_0}(x)} R_\theta, Df_{f^{j_0+1}(x)}, \ldots, Df_{f^m(x)} \}
$$
is a $(k,\UU)$-realizable sequence of length $m$ at $x$.
The product $L_{m-1} \cdots L_0$
sends $E^u(x)$ to $E^s(f^m(x))$, as required.

\subsubsection*{Second case}
We assume the following condition:
\begin{equation}\tag{II}
\text{there exist $j_0, j_1 \in \{0,1,\ldots,m-1\}$, $j_0<j_1$, s.t. }
\Delta (f^{j_0}(x), j_1-j_0) > c.
\end{equation}

Since $c>C^2$, we have $j_1 - j_0 > 1$.
Let $A = Df^{j_1-j_0}_{f^{j_0}x}$, and take unit vectors 
$s \in E^s(f^{j_0}(x))$ and $u \in E^u(f^{j_0}(x))$.
By (II), we have $\|A s\| / \| Au \| > c$.
Therefore lemma~\ref{linear1} gives a vector $\xi \in T_{f^{j_0}(x)}M$ such that
$$
|\theta_0| = \ang(\xi, E^s(f^{j_0}x)) \leq \alpha_2
\quad \text{and} \quad
|\theta_1| = \ang(Df^{j_1-j_0}_{f^{j_0}x}(\xi), E^u(f^{j_1}x)) \leq \alpha_2.
$$
The signs of $\theta_0$ and $\theta_1$ are chosen so that
$R_{\theta_0} ( E^s(f^{j_0} x) ) \ni \xi$ and
$R_{\theta_1} Df^{j_1-j_0}_{f^{j_0} x} (\xi) \in E^u (f^{j_1} x)$.
Applying lemma~\ref{tipo2}, we conclude that the sequence
$$
\{Df_{f^{j_0} x} R_{\theta_0},Df_{f^{j_0+1} x}, \dots ,
Df_{f^{j_1-2} x}, R_{\theta_1} Df_{f^{j_1-1} x} \},
$$
of length $j_1 - j_0$ at $f^{j_0}(x)$, is $(k_0,\UU)$-realizable.
Now define the sequence $\{L_0$, \ldots, $L_{m-1} \}$ of length $m$ at $x$
putting 
$L_{j_0} = Df_{f^{j_0} x} R_{\theta_0}$, 
$L_{j_1 - 1} = R_{\theta_1} Df_{f^{j_1} x}$ and all the others
$L_j = Df_{f^j x}$. 
By lemma~\ref{basico}, items 1 and 2, this is 
is a $(k,\UU)$-realizable sequence.
If the signs are appropriately chosen, then we have
$L_{m-1} \cdots L_0 (E^u(x)) = E^s(f^m(x))$.

\subsubsection*{Third case}
We suppose that we are not in the previous cases, that is 
we assume:
\begin{equation}\tag{not I}
\text{for every $j \in \{0,1,\ldots,m-1\}$,}\quad
\ang(E^u(f^j(x)), E^s(f^j(x)) \geq \alpha_1.
\end{equation}
and
\begin{equation}\tag{not II}
\text{for every $i,j \in \{0,1,\ldots,m-1\}$ with $i<j$,} \quad
\Delta (f^i(x), j-i) \leq c.
\end{equation}
We now use the assumption $x \in \Gamma_m(f)$, that is
$\Delta (x,m) \geq 1/2$.

\begin{claim}
For every $i,j \in \{0,1,\ldots,m-1\}$ with $i<j$,
$$
1/(2c^2) \leq \Delta (f^i(x), j-i) \leq c.
$$
\end{claim}

\begin{claim proof}
The second inequality is just (not~II). For the first:
$$
\Delta (f^i(x), j-i) = 
\Delta (f^j(x), m-j)^{-1} \cdot \Delta (x,m) \cdot \Delta (x,i)^{-1} \geq
1/(2c^2).
$$
\end{claim proof}

It follows from the claim, condition (not~I) and lemma~\ref{linear2} that
$$
\text{for every $j \in \{0,1,\ldots,m\}$,} \quad
\| Df_x^j \| \leq E.
$$
Choose numbers $\theta_0,\ldots,\theta_{m-1}$ such that 
$|\theta_j |\leq \beta$ and
$\sum \theta_j = \ang (E^u(x),E^s(x))$.
Define linear maps $L_j:T_{f^j x}M \to T_{f^{j+1} x}M$ by 
$$
L_j= Df^{j+1}_x R_{\theta_j} (Df^j_x)^{-1} \, .
$$
Then we have:
\begin{itemize}
\item[(i)] 
$L_{j-1}\ldots L_{0} = Df^j_x R_{\theta_{j-1}+\cdots+\theta_0}$, therefore
$L_{j-1}\ldots L_{0}(\mathbb{D})=Df_{x}^{j}(\mathbb{D})$, where
$\mathbb{D}$ is the unit disk in $T_x M$,
\item[(ii)]
$\| L_j - Df_{f^{j}x} \| \leq
\|Df_{f^{j}x} \| \left\| Df^j_x R_{\theta_j} (Df^j_x)^{-1} -I \right\| \leq
C \left\| Df^j_x \right\|^2 \| R_{\theta_j} - I \| < \varepsilon$. 
\end{itemize}
So we have constructed a sequence 
$\{L_0,\cdots,L_{m-1}\}$ of length $m$ at $x$ that, by lemma~\ref{tipo3}, 
is $(k, \UU)$-realizable. Furthermore, we have
$$
L_{m-1}\cdots L_{0}(E^u(x))=
Df^m_x R_{\sum \theta_j}(E^u(x))=
Df^m_x (E^s(x))=E^s(f^m x) \, ,
$$
as required.

\end{proof}

%%%%%%%%%%%%%%%%%%%%%%%%%%%%%%%%%%%%%%%%%%%%%%%%%%%%%%%%%%%%%%%%%%%%%%%%%
\subsection{Realizable sequences with small products} \label{s:prova ML}
%%%%%%%%%%%%%%%%%%%%%%%%%%%%%%%%%%%%%%%%%%%%%%%%%%%%%%%%%%%%%%%%%%%%%%%%%

In section \ref{s:Eu em Es} we have defined sets $\Gamma_m(f)$,
for $f\in\Diff$ and $m\in \NN$.
Define also the following $f$-invariant sets
\[
\Omega_m(f) = \bigcup_{n\in\ZZ}f^n(\Gamma_m(f)).
\]

\begin{lemma} \label{hiper bobo}
For $f\in \Diff$ and $m\in \NN$, let $H_m = \mathcal{O}^{+}(f) - \Omega_m(f)$.
If $H_m$ is not empty then its closure $\overline{H_m}$ is 
a hyperbolic set.
\end{lemma}

\begin{proof}
The proof is quite standard.
If $x\in H_m$ then $\Delta(x,m) \leq 1/2$ and $\Delta(x,mi) \leq 1/2^i$.
Since $\|Df\|$ is bounded, there are constants $K>1$ and $0<\tau<1$ such that
$$
\Delta(x,n) \leq K \tau^n \quad 
\text{for every $x\in H_m$ and $n\geq1$.}
$$
Fix $x\in H_m$ and let $v^u \in E^u(x)$, $v^s \in E^s(x)$ be unit vectors.
Then
\begin{multline*}
\| v^u - v^s \| \geq
\|Df^m_x\|^{-1} \| Df^m_x v^u - Df^m_x v^s \| \geq \\
\|Df^m_x\|^{-1} \big[ \| Df^m_x v^u \| - \| Df^m_x v^s \| \big] \geq
\|Df^m_x\|^{-1} \| Df^m_x v^u \| / 2.
\end{multline*}
This shows that there exists a constant $a>0$ such that 
$\theta(x) = \ang (E^u(x), E^s(x)) \geq a$ for every $x\in H_m$.

Let $x\in H_m$ and $n\geq 1$. Since $Df_x^n$ preserves the area form,
we have
\[
\sin \theta(x) =\left\| Df_{x}^{m}\mid _{E^{s}(x)}\right\| .\left\|
Df_{x}^{m}\mid _{E^{u}(x)}\right\| \sin \theta(f^n x). 
\]
Therefore 
\begin{alignat*}{2}
\left\| Df_x^n \mid _{E^s (x)} \right\| &= 
\sqrt{\frac{\sin \theta(x)}{\sin \theta(f^n x)} \Delta (x,n)} &\leq
\sqrt{\frac{K\tau^n}{\sin a}} =
K_{1}\tau _{1}^{n}, \\
\left\| Df_{x}^{-n}\mid _{E^u(x)} \right\| &= 
\sqrt{\frac{\sin \theta(x)}{\sin \theta(f^{-n}(x))} \Delta (f^{-n}(x),n)} \, &\leq
\sqrt{\frac{K\tau^n}{\sin a}} =
K_{1}\tau _{1}^{n}, 
\end{alignat*}
for some constants $K_{1}>0$, $0<\tau_{1}<1$. These inequalities imply
that the bundles $E^{s}$, $E^{u}$ are continuous on $H_m$, and have
a unique continuous extension to the closure.
\end{proof}

We will need the following result,
which may be thought as a quantitative Poincar\'{e}'s recurrence theorem.

\begin{lemma} \label{poincare}
Let $f\in \Diff$. Let $\Gamma \subset M $ be a measurable set
with $\mu (\Gamma )>0$ and let 
\[
\Omega =\bigcup_{n\in \ZZ }f^{n}(\Gamma )\, .
\]
Take $\gamma >0$. Then there exists a measurable function $N_{0}:\Omega
\rightarrow \NN$ such that for a.e. $x\in \Omega $, every $n\geq N_{0}(x)$
and every $t\in [0,1]$ there is some $\ell \in \{0,1,\ldots n\}$ such
that $f^{\ell}(x)\in \Gamma $ and $\left| \frac{\ell}{n}-t\right| <\gamma$.
\end{lemma}

\begin{proof}
Let $\chi _{\Gamma }$ be the characteristic function
of the set $\Gamma $. Consider the Birkhoff sums 
\[
s_{n}(x)=\sum_{j=0}^{n-1}\chi _{\Gamma }(f^{j}(x)) 
\]

\begin{claim}
For a.e. $x\in \Omega $, the limit
$\lim_{n\rightarrow \infty }\frac{1}{n}s_{n}(x)$ exists and is
positive.
\end{claim}

\begin{claim proof}
Birkhoff's theorem gives the existence; we are
left to show the positivity. Let $Z\subset \Omega $ be the set where the
limit is zero. Let $Z_{0}=Z\cap \Gamma $. The a.e.-defined 
$f$-invariant function 
\[
\varphi =\lim_{n\rightarrow \infty }\frac{1}{n}\sum_{j=0}^{n-1}
\chi_{Z_{0}}\circ f^{j} \leq \lim_{n\rightarrow \infty }\frac{s_{n}}{n}
\]
vanishes in $Z_{0}$. This means that the set $P=\{\varphi >0\}$ is disjoint
from $Z_{0}$.
On the other hand, $\varphi(x)>0$ implies that some iterate of $x$
is in $Z_{0}$.
Since $P$ is invariant, it follows that $P=\varnothing $. Thus
$\varphi =0$ a.e. and $\mu (Z_{0})=\int \varphi d\mu =0$. But $\Omega
=\bigcup_{n\in \ZZ }f^{n}(\Gamma )$ means that
$Z=\bigcup_{n\in \ZZ }f^{n}(Z_{0})$
and therefore $\mu (Z)=0$. The claim is proved.
\end{claim proof}

Take $x\in \Omega$. Let
$a=\lim\limits_{n\rightarrow \infty }\frac{1}{n}s_{n}(x)$.
Take $0<\varepsilon <a$ such that 
$\frac{a+\varepsilon }{a-\varepsilon }<1+\gamma/2$.
Choose (measurably) $n_{0}$ such that $n\geq
n_{0}\Rightarrow \left| \frac{s_{n}}{n}-a\right| <\varepsilon $. Finally,
take an integer
$$N_{0}(x)>\max \left\{ \frac{2n_{0}}{\gamma (a-\varepsilon )},\,
\frac{4}{\gamma }\right\} .$$

Now, by contradiction, suppose that for some $n\geq N_0(x)$ there exists $t\in[0,1]$
such that $f^{\ell}(x)\notin \Gamma $ for every $\ell \in (n(t-\gamma), n(t+\gamma))$.
Let $[\ell_1,\ell_2]$ be the maximal closed subinterval of 
$(n(t-\gamma), n(t+\gamma)) \cap [0,n]$ with integer endpoints.
Then $\ell_{2}-\ell_{1} > n\gamma - 2 > n\gamma /2$.
If $\ell_{1}\geq n_{0}$ then 
\[
a-\varepsilon <
\frac{s_{\ell_{2}}}{\ell_{2}} = 
\frac{s_{\ell_{1}}}{\ell_{2}} \leq 
\frac{s_{\ell_{1}}}{\ell_{1} + n\gamma/2 } \leq
\frac{s_{\ell_{1}}}{\ell_{1}(1+\gamma/2 )} <
\frac{a+\varepsilon }{1+\gamma/2 } < 
a-\varepsilon , 
\]
a contradiction. Therefore $\ell_{1}<n_{0}$.
We have $\ell_{2} > n \gamma/2 >\frac{n_{0}}{a-\varepsilon }>n_{0}$, thus 
\[
a-\varepsilon <
\frac{s_{\ell_{2}}}{\ell_{2}} <
\frac{\ell_{1}}{\ell_{1} + n\gamma/2 } <
\frac{n_{0}}{n\gamma/2} <
a-\varepsilon , 
\]
again a contradiction.
\end{proof}

In lemma \ref{Eu em Es} we have constructed realizable sequences that send $E^u$ in $E^s$.
Using this, we now construct realizable sequences whose products have ``small'' norms.

\begin{lemma} \label{ML simplificado}
Let $f\in \Diff$ be aperiodic and such that every 
hyperbolic set has zero measure, $\UU = \UU(f,\eps_0)$, $\delta >0$ and $0<k<1$.
Then there exists a measurable integer
function $N:M\rightarrow \NN$ such that 
for a.e. $x\in M$ and every integer $n\geq N(x)$ there exists 
a realizable sequence 
$$
\{ L_0,\ldots, L_{n-1} \} = \{ L_0^{(x,n)},\ldots, L_{n-1}^{(x,n)} \}
$$
of length $n$ at $x$ such that
$$
\| L_{n-1} \cdots L_0 \| < e^{\frac{4}{5}n\delta}.
$$
\end{lemma}

\begin{proof}
Let $k_0 \in (0, k)$.
Let $m\in \NN$, depending on $f$, $\eps_0$ and $k_0$ (in the place of $k$),
be given by lemma~\ref{Eu em Es}.
By lemma~\ref{hiper bobo}, the disjoint union
$\mathcal{O}^{0}(f) \sqcup \Omega_m(f)$ has full measure.
We will define the function $N: M \to \NN$ separately on $\mathcal{O}^{0}(f)$ and
$\Omega_m(f)$.

For each $x\in\mathcal{O}^{0}(f)$, take $N(x) \in \NN$ such that
$\| Df_x^n \| < e^{n\delta}$ for every $n\geq N(x)$.
Fixed some $n\geq N(x)$, we define
$L_j = L_j^{(x,n)} = Df_{f^j x}$ for $0 \leq j \leq n-1$.
By lemma~\ref{basico}, item 1, $\{L_0, \ldots L_{n-1} \}$ is a 
$(k,\UU)$-realizable sequence of length $n$ at $x$.

If $\mu(\mathcal{O}^{0}(f)) = 1$ then we are done.
Suppose from now on that $\mu(\Omega_m(f)) >0$.
Let $C > \log \sup_{g \in \UU, \, x\in M} \| Dg_x \|$.
Apply lemma~\ref{poincare} with
$\Gamma = \Gamma_m(f)$, $\Omega = \Omega_m(f)$,
and $\gamma =\frac{\delta }{20C}$ to find
$N_{0}(x)$, depending measurably on $x \in \Omega_m(f) $,
such that for every $n \geq N_{0}(x)$ and $t\in [0,1]$ there   
is $\ell \in \NN$ with $f^\ell(x) \in \Gamma_m(f)$ and  
$\left| \frac{\ell}{n}- t \right| <\gamma $.

Fix $x\in \Omega_m(f)$.
Consider the Lyapunov exponent
$\lambda =\lambda ^{+}(f,x) > 0$.
If $\lambda < \delta$ then it suffices to take
$N(x)$ large enough and defining $\{ L_j^{(x,n)} \}$ as the trivial realizable sequence
(as we did when defining $N$ on $\mathcal{O}^{0}(f)$).
Thus, we can assume that $\lambda \geq \delta $. 
If $x$ is contained in a certain full measure subset of $\Omega_m(f)$ 
(we will also assume this) 
then we can find $N_{1}(x)\geq N_{0}(x)$ such that if 
$j>\frac{1}{10}N_{1}(x)$ then 
\[
\left| \frac{1}{j}\log \left\| Df_{x}^{j}\mid _{E^{u}}\right\| -\lambda
\right| <\frac{\delta }{20}, 
\]
\[
\left| \frac{1}{j}\log \left\| Df_{x}^{j}\mid _{E^{s}}\right\| +\lambda
\right| <\frac{\delta }{20}, 
\]
\[
\left| \frac{1}{j}\log \sin \measuredangle
(E^{u}(f^{j}x),E^{s}(f^{j}x))\right| <\frac{\delta }{5}. 
\]

For each $j\geq 0,$ take unitary vectors $v_{j}^{u}\in E^{u}(f^{j}x)$ and $%
v_{j}^{s}\in E^{s}(f^{j}x)$. Let $\BB _{j}=\{v_{j}^{u},v_{j}^{s}\}$
be a basis of the space $T_{f^{j}x}M$. Denote also $\theta
_{j}=\measuredangle (v_{j}^{u},v_{j}^{s})$. If $A:T_{f^{j_{0}}x}M\rightarrow
T_{f^{j_{1}}x}M$ is any linear transformation, we will indicate by $\overline{%
A}$ its matrix relative to the bases $\BB _{j_{0}}$ and $\BB %
_{j_{1}}$. If 
\[
\overline{A}=\left( 
\begin{array}{cc}
a_{uu} & a_{su} \\ 
a_{us} & a_{ss}
\end{array}
\right) 
\]
is such a matrix, we will write $\left\| A\right\| _{\max }=\max \{\left|
a_{uu}\right| ,\left| a_{su}\right| ,\left| a_{us}\right| ,\left|
a_{ss}\right| \}$. We claim that there is a constant $K>1$ such that 
\begin{eqnarray*}
\left\| A\right\| &\leq &K(\sin \theta _{j_{0}})^{-1}\left\| A\right\|
_{\max }, \\
\left\| A\right\| _{\max } &\leq &K(\sin \theta _{j_{1}})^{-1}\left\|
A\right\|
\end{eqnarray*}
for every such $A$. To prove this fact, look at the matrix of the change of
bases from $\BB _{j}$ to the orthonormal basis $%
\{v_{j}^{u},(v_{j}^{u})^{\bot }\}$, that is, 
\[
\left( 
\begin{array}{cc}
1 & \cos \theta _{j} \\ 
0 & \sin \theta _{j}
\end{array}
\right) . 
\]
The norm of its inverse is of the order of $(\sin \theta _{j})^{-1}$.

Take $N(x) \geq N_{1}(x)$ such that
$$
\frac{m}{N(x)}<\gamma
\quad \text{and} \quad
\frac{K^{2} e^{Cm}}{\sin \theta_0} < e^{\frac {\delta }{5}N(x)}\, .
$$ 
(Notice $\theta_0$ is a measurable function of $x$.)
This defines the measurable function $N$ on the set $\Omega_m(f)$.
Now fix some $n\geq N(x)$. 
Let $\ell$ be as given by lemma~\ref{poincare} with
$t=\frac{1}{2}$, that is, such that 
$f^\ell(x) \in \Gamma_m(f)$ and  
$\left| \frac{\ell}{n}-\frac{1}{2} \right| <\gamma $.
By lemma~\ref{Eu em Es}, there is a 
$(k_0,\UU)$-realizable sequence
$\{L'_0, \ldots L'_{m-1} \}$ at $f^\ell(x)$ of length $m$
such that the product $\mathcal{L}' =  L'_{m-1} \cdots L'_0$
satisfies $\mathcal{L}' (E^u(f^\ell (x))) = E^s (f^{\ell + m}(x))$.
We now define the sequence $\{L_0, \ldots, L_{n-1}\}$ of 
length $n$ at $x$ as equal to 
$$
\{ Df_x, \ldots, Df_{f^{\ell-1}x}, L'_0, \ldots L'_{m-1} , 
Df_{f^{\ell +m + 1} x}, \ldots, Df_{f^n x} \}
$$
By lemma~\ref{basico}, it is a $(k,\UU)$-realizable sequence.
To complete the proof, we must estimate the norm of
$$
\mathcal{L} = L_{n-1} \cdots L_0 = 
Df_{f^{\ell + m}(x)}^{n-\ell-m} \cdot \mathcal{L}' \cdot Df_{x}^{\ell}.
$$

Write 
\[
\overline{Df_{x}^{\ell}}=\left( 
\begin{array}{cc}
e^{\mu _{1}n} & 0 \\ 
0 & e^{-\mu _{2}n}
\end{array}
\right) ,\quad 
\overline{Df_{f^{\ell + m}(x)}^{n-\ell-m}}=\left( 
\begin{array}{cc}
e^{\mu _{3}n} & 0 \\ 
0 & e^{-\mu _{4}n}
\end{array}
\right) . 
\]

\begin{claim}
For $i=1,2,3,4$ we have
$\left| \mu _{i}-\frac{\lambda }{2}\right| <\frac{\delta }{5}.$
\end{claim}

\begin{claim proof}
We have
$\mu _{1}=\frac{1}{n}\log \left\| Df_{x}^{\ell}\mid _{E^{u}}\right\|$. But
$$
\left| \frac{1}{\ell}\log \left\| Df_{x}^{\ell}\mid _{E^{u}}\right\| -\lambda
\right| <\frac{\delta }{20}
\quad \text{and} \quad
\left| \frac{\ell}{n}-\frac{1}{2} \right| <\gamma \, .
$$
Also, $\lambda \leq C$. Through $| ab-a'b'| \leq |a-a'| |b| +|b-b'| |a'|$,
we obtain
$$ 
\left| \mu _{1}-\frac{\lambda }{2}\right| <
\frac{\delta }{20}+C\gamma =
\frac{2 \delta }{20}\, .
$$
Similarly for $\mu _{2}$. Now,
$$
\mu _{3}=
\frac{1}{n}\log \left\| Df_{f^{\ell+m}(x)}^{n-\ell-m}\mid_{E^{u}}\right\| =
\frac{1}{n}\log \left\| Df_{x}^{n}\mid_{E^{u}}\right\| -
\frac{1}{n}\log \left\| Df_{x}^{\ell+m}\mid _{E^{u}}\right\| \, .
$$
Since 
$$
\left|\frac{1}{\ell+m}\log \left\| Df_{x}^{\ell+m}\mid _{E^{u}}\right\| 
- \lambda \right| <\frac{\delta }{20}
\quad \text{and} \quad
\left| \frac{\ell+m}{n}-\frac{1}{2} \right| <2\gamma \, ,
$$
we have
$$
\left|\frac{1}{n}\log \left\| Df_{x}^{\ell+m}\mid _{E^{u}}\right\| 
- \frac{\lambda }{2} \right| < 
\frac{\delta }{20}+2C\gamma =
\frac{3 \delta }{20} \, .
$$
Thus
$$
\left| \mu _{3}-\frac{\lambda }{2}\right| \leq
\left| \frac{1}{n}\log \left\| Df_{x}^{n}\mid_{E^{u}}\right\| 
- \lambda \right| +
\left| \frac{1}{n}\log \left\| Df_{x}^{\ell+m}\mid _{E^{u}}\right\| 
- \frac{\lambda }{2} \right| <
\frac{\delta }{20} + \frac{3 \delta }{20} =
\frac{\delta }{5} \, .
$$
Similarly for $\mu _{4}$.
\end{claim proof}

Since the sequence $\{L'_0, \ldots, L'_{m-1} \}$ is $(k,\UU)$-realizable,
we have $\| L'_i \| < e^C$ for each $i$.
In particular, $\| \mathcal{L}' \| < e^{Cm}$ and
$$
\| \mathcal{L}' \| _{\max }\leq
K(\sin \theta_{\ell+m})^{-1} \| \mathcal{L}' \| <
K e^{Cm} e^{\frac{1}{5}\delta n}\, ,
$$
since $\sin \theta_{\ell+m} > e^{-\frac15 \delta (\ell +m)} > e^{-\frac15 \delta n}$.
Now write 
\[
\overline{\mathcal{L}' }=\left( 
\begin{array}{cc}
b_{uu} & b_{su} \\ 
b_{us} & b_{ss}
\end{array}
\right) . 
\]
Since $\mathcal{L}'(E^u_{f^{\ell}(x)})=E^s_{f^{\ell+m}(x)}$, we have $b_{uu}=0$.
Computing the matrix product
$\overline{\mathcal{L}} = 
\overline{Df_{f^{\ell + m}(x)}^{n-\ell-m}} \cdot 
\overline{\mathcal{L}'} \cdot 
\overline{Df_{x}^{\ell}}$,
we obtain
\[
\overline{\mathcal{L}}=\left( 
\begin{array}{cc}
0 & e^{(-\mu _{2}+\mu _{3})n}b_{su} \\ 
e^{(\mu _{1}-\mu _{4})n}b_{us} & e^{(-\mu _{2}-\mu _{4})n}b_{ss}
\end{array}
\right) . 
\]
We have $\max \{ \mu_1-\mu_4, \; -\mu_2+\mu_3, \; -\mu_2-\mu_4 \} < \frac{2}{5}\delta$.
Therefore
\[
\| \mathcal{L} \| _{\max }\leq 
e^{\frac{2}{5}\delta n} \| \mathcal{L}' \| _{\max }\leq Ke^{Cm} e^{\frac{3}{5}\delta n}
\]
and
\[
\| \mathcal{L} \| \leq 
K     (\sin \theta _{0})^{-1} \| \mathcal{L} \| _{\max }\leq 
K^{2} (\sin \theta _{0})^{-1} e^{Cm} e^{\frac{3}{5}\delta n}<
e^{\frac{4}{5}\delta n}\,.
\]
This concludes the proof.
\end{proof}

The Main Lemma follows easily from lemma~\ref{ML simplificado}:

\begin{proof}[Proof of the Main Lemma]
Of course, we can suppose that 
$\UU = \UU (f,\varepsilon _{0})$ for some $\varepsilon _{0}>0$.
Applying lemma~\ref{ML simplificado}, we find the measurable function
$N: M \to \NN$.
Fix $x$ (in a full measure set) and $n \geq N(x)$.
Then lemma~\ref{ML simplificado}
also gives us a $(k,\UU)$-realizable sequence
$\{L_0, \ldots, L_{n-1} \}$ of length $n$ at $x$ such that
\begin{equation}\tag{$\ast$}
\| L_{n-1} \cdots L_0 \| < e^{\frac{4}{5}n\delta}.
\end{equation}
Take $\gamma >0$ very small (depending on $n$).
By the definition of a realizable sequence,
there exists $r=r(x,n)$ with the following property:
For every disk $U = B_{r'}(x)$, with $0<r'<r$ there exist
$g\in \UU(f,\eps_0)$ and a compact set $K\subset U$ such that:
\begin{itemize}
\item[(i)] $g$ equals $f$ outside the set
$\bigsqcup_{j=0}^{n-1}f^{j}(\overline{U})$ and the iterates 
$f^{j}(\overline{B_{r}(x)}),\, 0 \leq j \leq n-1$ are two-by-two
disjoint,

\item[(ii)] $\frac{\mu (K)}{\mu (U)}>k$,

\item[(iii')] if $y\in K$ then $\left\| Dg_{g^{j}y} - L_{j}\right\| <\gamma 
$ for every $j = 0,1,\ldots,n-1$.

\end{itemize}
If $\gamma$ is small enough then ($\ast$) and (iii') imply that
%\begin{equation} \tag{iii}
%\text{if $y\in K$ then } \| Dg^n_y \| < e^{n\delta}.
%\end{equation}
\begin{itemize}
\item[(iii)] if $y\in K$ then $\| Dg^n_y \| < e^{n\delta}$.
\end{itemize}
This completes the proof of the Main Lemma.
\end{proof}

%%%%%%%%%%%%%%%%%%%%%%%%%%%%%%%%%%%%%%%%%%%%%%%%%%%%%%%%%%%%%%%%%%
%%%%%%%%%%%%%%%%%%%%%%%%%%%%%%%%%%%%%%%%%%%%%%%%%%%%%%%%%%%%%%%%%%
\section{Proof of Theorem A} \label{s.prova teor A}
%%%%%%%%%%%%%%%%%%%%%%%%%%%%%%%%%%%%%%%%%%%%%%%%%%%%%%%%%%%%%%%%%%
%%%%%%%%%%%%%%%%%%%%%%%%%%%%%%%%%%%%%%%%%%%%%%%%%%%%%%%%%%%%%%%%%%

%%%%%%%%%%%%%%%%%%%%%%%%%%%%%%%%%%%%%
\subsection{Preliminary definitions}
%%%%%%%%%%%%%%%%%%%%%%%%%%%%%%%%%%%%%

Fix some measure preserving diffeomorphism $f$. If a measurable set $%
A\subset M$ and $n\in \NN$ are such that the sets $A,f(A),\ldots
,f^{n-1}(A)$ are disjoint then we call the set 
\[
T=A\cup f(A)\cup \cdots \cup f^{n-1}(A) 
\]
a \emph{tower} for $f$. The number $n$ is called the \emph{height} of
the tower and the set $A$ is called its \emph{base}. A \emph{castle} is
a finite or countable union of two-by-two disjoint towers. The 
\emph{base of the castle} is the union of the bases of its towers.

Given $f$ and a positive measure set $A$ $\subset M$, consider the (a.e.
finite) return time $\tau :A\rightarrow \NN$ defined by $\tau (x)=\inf
\{n\geq 1;\;f^{n}(x)\in A\}$. If we denote $A_{n}=\tau ^{-1}(n)$ then $%
T_{n}=A_{n}\cup f(A_{n})\cup \cdots \cup f^{n-1}(A_{n})$ is a tower. Let $%
Q=\bigcup_{n\in \ZZ }f^{n}(A)$. Then $Q$ is $f$-invariant and it is a
castle with base $A$ and towers $T_{n}$. We will call $Q$ the 
\emph{Kakutani castle} with base $A$.

We will need also the following:

\begin{lemma}
For every aperiodic invertible measure preserving transformation $f$ on a
probability space $X$, every subset $U\subset X$ of positive measure and
every $n\in \NN$, there exists a positive measure set $V\subset U$
such that $V,f(V),\ldots ,f^{n}(V)$ are two-by-two disjoint.
Besides, $V$ can be chosen maximal on ``the measure-theoretical sense''.
(This means that no set that includes $V$ and has larger measure than $V$
has the stated properties.)
\end{lemma}

\begin{proof}
We follow \cite[page 70]{Halmos}.
Take $U_{1}\subset U$ such that $\mu (U_{1}\vartriangle f(U_{1}))>0$ (it exists
because otherwise a.e. point of $U$ would be fixed). Then $%
V_{1}=U_{1}-f(U_{1})$ has positive measure and $V_{1}\cap
f(V_{1})=\varnothing $. Take $U_{2}\subset V_{1}$ such that $\mu
(U_{2}\vartriangle f^{2}(U_{2}))>0$ and let $V_{2}=U_{2}-f^{2}(U_{2})$.
Continuing in this way we will find $V=V_{n}$ such that $V,f(V),\ldots
,f^{n}(V)$ are two-by-two disjoint. Suppose that the set $R_{V}=U-%
\bigcup_{j=-n}^{n}f^{j}(V)$ has positive measure; otherwise $V$ is maximal.
Take $V^{\prime }\subset R_{V}$ such that $V^{\prime },f(V^{\prime }),\ldots
,f^{n}(V^{\prime })$ are two-by-two disjoint. Continue in this way by
transfinite induction. Since a disjoint class of positive measure sets is
countable, the process will terminate at some countable ordinal. Hence we
find a measurable set $V\subset U$ such that $\mu (R_{V})=0$.
\end{proof}

Now we will prove proposition \ref{teorema} and thus theorem A.

\subsection{First step. Construction of a castle $Q$}

Let $f\in \Diff$ be
aperiodic and such that every hyperbolic set for $f$ has zero measure.
Let $\UU $ be a neighborhood of $f$ in $\Diff$ and let $\delta >0$.

Take $0<k<1$ such that $1-k<\delta ^{2}$. Apply the Main Lemma to get a
measurable function $N:M\rightarrow \NN$ with the properties stated there.
We define the sets $P_{n}=\left\{ x\in M;\;N(x)\leq n\right\} $
for $n\in \NN$. Obviously, $\mu (P_{n})\rightarrow 1$.
Fix $H\in \NN$ such that $\mu
(P_{H}^{C})<\delta ^{2}$, where $P_{H}^{C}$ denotes the complementary set
$M-P_{H}$.

Take $B\subset P_{H}$ such that $B,f(B),\ldots ,f^{H-1}(B)$ are two-by-two
disjoint and such that $B$ is maximal in measure-theoretical sense. Consider
the following $f$-invariant set:
\[
\widehat{Q}=\bigcup_{n\in \ZZ }f^{n}(B). 
\]
$\widehat{Q}$ is the Kakutani castle with base $B$. Notice that $\widehat{Q}%
\supset P_{H}$ mod 0 (by maximality of $B$) and hence $\mu (\widehat{Q}%
^{C})<\delta ^{2}$. Let $Q\subset \widehat{Q}$ be the (finite) castle
consisting of all the towers of $\widehat{Q}$ with heights at most $3H$
floors. The following property will be important later:

\begin{lemma} \label{3delta2}
$\mu (\widehat{Q}-Q)\leq 3\mu (P_{H}^{C})<3\delta ^{2}$.
\end{lemma}

\begin{proof}
Write the castle as $\widehat{Q}%
=\bigsqcup_{i=H}^{\infty }T_{i}$ where $B=\bigsqcup_{i=H}^{\infty }B_{i}$ is
the base and $T_{i}=\bigsqcup_{j=0}^{i-1}f^{j}(B_{i})$ is the tower of
height $i$ floors. Take $i\geq 2H$ and $H\leq j\leq i-H$. The sets $%
f^{j}(B_{i}),\ldots f^{j+H-1}(B_{i})$ are disjoint and do not intersect $%
B\sqcup \cdots \sqcup f^{H-1}(B)$. Since $B$ is maximal, we conclude that (see
figure 1)
\[
i\geq 2H,\ H\leq j\leq i-H \Longrightarrow f^{j}(B_{i})\subset P_{H}^{C}\mbox{
mod 0} 
\]
(otherwise we could replace $B$ by $B \sqcup (f^{j}(B_{i})\cap P_{H})$,
contradicting the maximality of $B$.)
\begin{figure}[hbt]
\begin{center}
\psfrag{H}{$H$}
\psfrag{3H}{$3H$}
\psfrag{castle Q}{castle $Q$}
\psfrag{base B}{base $B$}
\psfrag{this part is}{this part is}
\psfrag{contained in PHC}{contained in $P_H^C$}
\includegraphics[height=2.0in]{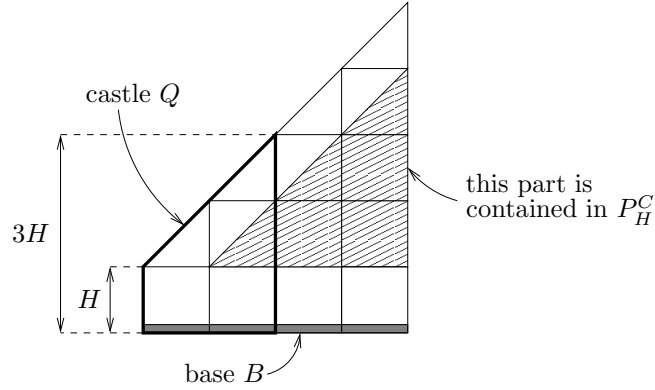}
\caption{The castle $\widehat{Q}$.}
\end{center}
\end{figure}
Thus
\[
i\geq 2H \Rightarrow \mu (T_{i}\cap P_{H}^{C})\geq
\sum_{j=H}^{i-H} \mu(f^j(B_i)) = \frac{i-2H+1}{i}\mu (T_{i}). 
\]
In particular, 
\[
i\geq 3H+1 \Rightarrow \mu (T_{i}\cap P_{H}^{C})>\frac{1}{3}\mu (T_{i}) 
\]
and so
\[
\mu (\widehat{Q}-Q)=
\sum_{i=3H+1}^{\infty }\mu (T_{i})\leq
\sum_{i=3H+1}^{\infty }3\mu (T_{i}\cap P_{H}^{C})=
3\mu \left(P_H^C \cap \bigsqcup_{i=3H+1}^{\infty }T_{i}\right) \leq
3\mu (P_{H}^{C}).
\]
\end{proof}

%%%%%%%%%%%%%%%%%%%%%%%%%%%%%%%%%%%%%%%%%%%%%%%%%%%%%%%%%%%%%%%%%%%%%%%%%%%
\subsection{Second step. Construction of the perturbed diffeomorphism $g$}
%%%%%%%%%%%%%%%%%%%%%%%%%%%%%%%%%%%%%%%%%%%%%%%%%%%%%%%%%%%%%%%%%%%%%%%%%%%

Let $0<\gamma <\delta ^{2}H^{-1}$.
By the regularity of the measure $\mu$, one can
find a compact castle $J\subset Q$ such that $\mu (Q-J)<\gamma $ and  of
the \emph{same type} as $Q$. 
Saying so we mean that the castles have the same number of towers and
that the towers have the same heights.
Since $J$ is compact, we can find an \emph{open} castle $V$ containing $J$
with $\mu (V-J)<\gamma $ and also with the same type as $J$.
Hence ($\vartriangle $ denotes symmetric difference) 
\[
\mu (V\vartriangle Q)=\mu (V-Q)+\mu (Q-V)\leq \mu (V-J)+\mu (Q-J)<2\gamma . 
\]

Denote by $S$ the base of the castle $V\cap Q$. For each $x\in S$, let 
$n(x)$ be the height of the tower that contains $x$.
We have $n(x)\geq H \geq N(x)$. Hence the Main Lemma gives, for a.e.
$x\in S$, a radius $r(x)=r(x,n(x))$.

Reducing $r(x)$ if needed, we can suppose that the disk 
$\overline{B}_{r(x)}(x)$ is contained in the base of a tower in $V$,
for a.e. $x\in S$.
Using Vitali's covering lemma, we can find a finite collection of disjoint
disks $U_{i}=B_{r_{i}}(x_{i})$ with $x_{i}\in S$ and $0<r_{i}<r(x_{i})$ such
that 
\begin{equation}
\frac{\mu \left( S-\bigsqcup \overline{U_{i}}\right) }{\mu \left( S\right) }%
<\gamma \, .
\tag{$\ast$}
\end{equation}
(Actually, Vitali's lemma allows us only to a.e.-cover the set $S$ restricted to each
chart domain, so we have to cover $S$ by chart domains first.)

Let $n_i = n(x_i)$. Notice that $n(x)=n_i$ for all $x\in U_i$.
By the Main Lemma, for each $i$ we can find a compact set $K_{i}\subset
U_{i} $ and $g_{i}\in \UU $ such that
\begin{itemize}
\item[(i)] $\frac{\mu (K_{i})}{\mu (U_{i})}>k$,

\item[(ii)] $g_{i}$ equals $f$ outside the set $%
\bigsqcup_{j=0}^{n_{i}-1}f^{j}(U_{i})$,

\item[(iii)] $x\in K_{i}\Rightarrow \big\| D(g_{i}^{n_i})_{x} \big\|
<e^{\delta n_i}.$
\end{itemize}

Let $g$ be equal to $g_{i}$ in the set 
$\bigsqcup_{j=0}^{n_{i}-1}f^{j}(U_{i})$,
for each $i$, and be equal to $f$ outside. 
Since those sets are disjoint, $g\in \Diff$ is a well-defined diffeomorphism
Moreover, $g\in \UU =\UU (f,\varepsilon _{0})$.

Since each $U_i$ is contained in the base of a tower in the castle $V$,
$V$ is also a castle for $g$.
Moreover, we can define a $g$-castle $U$ of the same type as $V$ with base
$\bigsqcup U_{i}$. Analogously, let $K$ be the $g$-castle of the same 
type as $V$ with base $\bigsqcup K_{i}$.
Then, by definitions, $K\subset U\subset V$.
We have $\frac{\mu (K)}{\mu (U)}>k$. Also, ($\ast$) implies
$\frac{\mu (U)}{\mu (V\cap Q)}>1-\gamma$. Thus 
\[
\mu (V-U)\leq \mu (V\cap Q-U)<\gamma \mu (V\cap Q)<\gamma 
\]
and 
\[
\mu (U\vartriangle Q)\leq \mu (U\vartriangle V)+\mu (V\vartriangle
Q)<3\gamma . 
\]

Summarizing, we have constructed a diffeomorphism $g\in \UU $ and $g$%
-castles $K\subset U$ of the same type as the castle $Q$ and such that
\begin{itemize}
\item[(i)] $\mu (U\vartriangle Q)<3\gamma $ and $\mu (U-K)<1-k$,

\item[(ii)] $g$ equals $f$ outside $U$,

\item[(iii)] if $x$ is in base of $K$ and $n(x)$ is the height of the
tower of $K$ that contains $x$ then we have 
$ \big\| Dg_{x}^{n(x)}\big\| <e^{\delta n(x)} $.
\end{itemize}

%%%%%%%%%%%%%%%%%%%%%%%%%%%%%%%%%%%%%%%%%%%%%%%%%%%%%%%%%%%%%%%%%%%%%%%%
\subsection{Third step. Estimation of $\mathrm{LE}(g)$}
%%%%%%%%%%%%%%%%%%%%%%%%%%%%%%%%%%%%%%%%%%%%%%%%%%%%%%%%%%%%%%%%%%%%%%%%

We claim that $\mathrm{LE}(g)$ is small. To show it, we will use the
following property from proposition \ref{semicont}: 
\[
g\in \Diff\Rightarrow \mathrm{LE}(g)=\inf_{N\geq 1}%
\frac{1}{N}\int_{M}\log \left\| Dg^{N}\right\| d\mu . 
\]
This allows one to conclude in finite time (i.e., without taking limits)
that the integrated Lyapunov exponent is small.

Set $N=\delta ^{-1}H$. (Of course, we can assume that $\delta ^{-1}\in 
\NN $.) We define the following ``good set'' 
\[
G=\bigcap_{j=0}^{N-1}g^{-j}(K). 
\]
We will show that this set has almost full measure:

\begin{lemma}
$\mu(G^C) < 16\delta$.
\end{lemma}

\begin{proof}
Consider the following sets. 
\[
S_{1}=U-K,\quad S_{2}=Q-U,\quad S_{3}=\widehat{Q}-Q,\quad S_{4}=M-\widehat{Q}%
. 
\]
Then $K^{C}\subset S_{1}\cup S_{2}\cup S_{3}\cup S_{4}$ and 
\[
G^{C}\subset \bigcup_{i=1}^{4}\bigcup_{j=0}^{N-1}g^{-j}(S_{i}). 
\]
We will estimate the measure of each set $\cup_{j=0}^{N-1}g^{-j}(S_{i})$
separately.

First, we have $\mu (S_{2}) \leq \mu (U\vartriangle Q)<3\gamma $, and so
\[
\mu \left( \bigcup_{j=0}^{N-1}g^{-j}(S_{2})\right) \leq N\mu
(S_{2})<3N\gamma <3\delta . 
\]

Before continuing, we point out that if $X\subset M$ is any measurable set then 
\begin{align*}
\mu \left( \bigcup_{j=0}^{N-1}g^{-j}(X)\right) &\leq 
\mu \left( X \cup \bigcup_{j=1}^{N-1} (g^{-j}(X)-g^{-(j-1)}(X)) \right) \\
&= \mu (X)+(N-1)\mu \left(g^{-1}(X)-X\right) . 
\end{align*}
Moreover,
$\mu \left( g^{-1}(X)-X\right) =\mu \left( g(X)-X\right)=\mu \left( X-g(X)\right)$,
simply because $\mu$ is $g$-invariant.

We have
\[
\mu \left( \bigcup_{j=0}^{N-1}g^{-j}(S_{1})\right) \leq 
\mu (S_{1}) + N \mu(S_1 - g(S_1)) \, .
\]
The set $S_{1}=U-K$ is a $g$-castle whose towers have heights at least $H$.
Hence its first floor, which contains the set $S_1 - g(S_1)$,
measures at most $\frac{1}{H}\mu (S_{1})$. Moreover, $\mu(S_1)< 1-k = \delta^2$.
Substituting these estimates, we get:
\[
\mu \left( \bigcup_{j=0}^{N-1}g^{-j}(S_{1})\right) <
\delta^2 + \frac{N\delta^2}{H} < \delta^2 + \delta < 2\delta\, .
\]

We are going to treat the case $i=3$ similarly. By lemma \ref{3delta2},
$\mu(S_3)<3\delta^2$.
The set $S_{3}=\widehat{Q}-Q$ is an $f$-castle whose towers have heights at
least $3H$. Hence $\mu \left( f(S_{3})-S_{3}\right) \leq \frac{1}{3H}\mu
(S_{3})$. Since $\ g$ and $f$ coincide in $U^{C}$, we have 
\[
g(S_{3})=g(S_{3}-U)\cup g(S_{3}\cap U)\subset f(S_{3})\cup g(S_{3}\cap U). 
\]
So
\[
\mu \left( g(S_{3})-S_{3}\right) \leq \mu \left( f(S_{3})-S_{3}\right) +\mu
\left( g(S_{3}\cap U)\right) \leq \frac{\mu (S_{3})}{3H}+\mu (U-Q)<\frac{%
\delta ^{2}}{H}+3\gamma 
\]
and 
\[
\mu \left( \bigcup_{j=0}^{N-1}g^{-j}(S_{3})\right) <3\delta ^{2}+N\left( 
\frac{\delta ^{2}}{H}+3\gamma \right) <3\delta ^{2}+\delta +3N\gamma
<7\delta . 
\]

The set $S_{4}=M-\widehat{Q}$ is $f$-invariant, thus 
\[
g(S_{4})-S_{4}=\left[ g(S_{4}-U)\cup g(S_{4}\cap U)\right] -S_{4}\subset
g(S_{4}\cap U)\subset g(U-Q) .
\]
Now, $\mu(S_4)< \delta^2$ (see the definition of $\widehat{Q}$) and
$\mu(U-Q) \leq \mu (U\vartriangle Q)<3\gamma$. Thus
\[
\mu \left( \bigcup_{j=0}^{N-1}g^{-j}(S_{4})\right) <\mu (S_{4})+N\mu
(U-Q)<\delta ^{2}+3N\gamma <4\delta . 
\]

Putting all estimates together, we get
$ \mu (G^{C})<2\delta +3\delta +7\delta +4\delta =16\delta$.
\end{proof}

Now we will show that $\frac{1}{N}\log \left\| Dg^{N}\right\| $ is small
inside the set $G$. Let $K_{0}$ be the base of the castle $K$. For $y\in
K_{0}$, let $n(y)$ be the height of the tower containing $y$. We know that 
\[
\big\| Dg_{y}^{n(y)}\big\| <e^{\delta n(y)}. 
\]
Now take $x\in G$. Since the heights of $K$-towers are less than $3H$, we
can write 
\[
N=j_{1}+n_{1}+n_{2}+\cdots +n_{i}+j_{2} 
\]
such that $0\leq j_{1},j_{2}<3H,$ $1\leq n_{1},\ldots ,n_{i}<3H$, and the
points 
\[
f^{j_{1}}(x),f^{j_{1}+n_{1}}(x),\ldots ,f^{j_{1}+n_{1}+\cdots +n_{i}}(x) 
\]
are exactly the points of the orbit segment 
$x, f(x),\ldots,f^{N-1}(x)$ which belong to $K_{0}$.
Hence, if $C>\sup_{y\in M}\left\| Dg_{y}\right\| $
then 
\[
\left\| Dg_{x}^{N}\right\| <C^{j_{1}}e^{\delta n_1}\cdots e^{\delta n_i}
C^{j_{2}}<C^{6H}e^{\delta N} 
\]
and so
\[
\frac{1}{N}\log \left\| Dg_{x}^{N}\right\| <\frac{6H\log C}{N}+\delta
=\left( 6\log C+1\right) \delta . 
\]
Therefore
\begin{eqnarray*}
\mathrm{LE}(g) 
&\leq &\int_{M}\frac{1}{N}\log \left\| Dg^{N}\right\| d\mu \\
&=&\int_{G}\frac{1}{N}\log \left\| Dg^{N}\right\| d\mu 
+\int_{G^{C}}\frac{1}{N}\log \left\| Dg^{N}\right\| d\mu \\
&<&\left( 6\log C+1\right) \delta \cdot \mu (G)+(\log C)\mu (G^{C}) \\
&<&(22\log C+1)\delta \, .
\end{eqnarray*}
This proves theorem A.

%%%%%%%%%%%%%%%%%%%%%%%%%%%%%%%%%%%%%%%%%%%%%%%%%%%%%%%%%%%%%%
%%%%%%%%%%%%%%%%%%%%%%%%%%%%%%%%%%%%%%%%%%%%%%%%%%%%%%%%%%%%%%
\section{Generic dichotomy for continuous cocycles}
%%%%%%%%%%%%%%%%%%%%%%%%%%%%%%%%%%%%%%%%%%%%%%%%%%%%%%%%%%%%%%
%%%%%%%%%%%%%%%%%%%%%%%%%%%%%%%%%%%%%%%%%%%%%%%%%%%%%%%%%%%%%%

%%%%%%%%%%%%%%%%%%%%%%%%%%%%%%%%%%%%
\subsection{Cocycles}
%%%%%%%%%%%%%%%%%%%%%%%%%%%%%%%%%%%%

Let $(X,\mu )$ be a non-atomic probability space and $T:X\hookleftarrow $ an
automorphism of it. Denote by $\sldoisr$ the group of two-by-two
real matrices with unit determinant. Let
\[
L^{\infty }(X,\sldoisr)=\left\{ A:X\rightarrow \sldoisr
\mbox{ measurable and essentially bounded}\right\} 
\]
and consider in this space the following metric 
\[
\left\| A-B\right\| _{\infty }=\mathrm{ess}\sup \left\| A(x)-B(x)\right\| . 
\]

Given $A\in L^\infty (X,\sldoisr)$, we denote for $x\in X$ and $n\in \ZZ$,
$$
A^n (x)=
\left\{
\begin{array}{l}
A(T^{n-1} x) \cdots A(x)                 \\ 
I                                        \\
\left[A(T^n x)\right]^{-1} \cdots \left[A(T^{-1}x)\right]^{-1} 
\end{array}
\right.
\begin{array}{l}
\text{if $n>0$,} \\ 
\text{if $n=0$,} \\
\text{if $n<0$.}
\end{array}
$$
Notice that the following relation, called the cocycle identity, 
\[
A^{n+m}(x)=A^{n}(T^{m}x) \cdot A^{m}(x) 
\]
is satisfied for $m,n\in \ZZ $. With some abuse, we call the function $A$
a \emph{cocycle}.

Oseledets' theorem may be summarized in the present case as follows:

\begin{theorem}
Let $T$ and $A$ be as above. Then there exists a measurable function $%
x\mapsto \lambda ^{+}(x)\geq 0$ such that 
\[
\lambda ^{+}(x)=\lim_{n\rightarrow +\infty }\frac{1}{n}\log \left\|
A^{n}(x)\right\| 
\]
for $\mu $-a.e. $x\in X$. Moreover, if $\mathcal{O}^{+}=\{x;\;\lambda
^{+}(x)>0\}$ has positive measure then for a.e. $x\in \mathcal{O}^{+}$ there
is a splitting $\RR ^{2}=E^{u}(x)\oplus E^{s}(x)$, depending measurably
on $x$, such that for $v\in \RR ^{2}-\{0\}$%
\begin{eqnarray*}
\lim_{n\rightarrow +\infty }\frac{1}{n}\log \left\| A^{n}(x) \cdot v\right\|
&=&\left\{ 
\begin{tabular}{l}
$\lambda ^{+}(x)$ if $v\notin E^{s}(x)$ \\ 
$-\lambda ^{+}(x)$ if $v\in E^{s}(x)$%
\end{tabular}
\right. , \\
\lim_{n\rightarrow -\infty }\frac{1}{n}\log \left\| A^{n}(x) \cdot v\right\|
&=&\left\{ 
\begin{tabular}{l}
$-\lambda ^{+}(x)$ if $v\notin E^{u}(x)$ \\ 
$\lambda ^{+}(x)$ if $v\in E^{u}(x)$%
\end{tabular}
\right. ,
\end{eqnarray*}
\[
\lim_{n\rightarrow \pm \infty }\frac{1}{n}\log \sin \measuredangle
(E^{u}(T^{n}x),E^{s}(T^{n}x))=0. 
\]
\end{theorem}

As before, we define the integrated Lyapunov exponent, 
\[
\mathrm{LE}(A)=\int_{X}\lambda ^{+}d\mu =\inf_{n\in \NN}\frac{1}{n}%
\int_{X}\log \left\| A^{n}\right\| d\mu . 
\]

We now define the notion of uniform hyperbolicity for cocycles.

\begin{definition}
A cocycle $A\in L^{\infty }(X,
\sldoisr )$ over $T:(X,\mu )\hookleftarrow $ is called \emph{uniformly
hyperbolic} if there exists, for a.e. $x\in M$, a splitting $\RR %
^{2}=E^{u}(x)\oplus E^{s}(x),$ which is measurable with respect to $x$, such
that:
\begin{itemize}
\item[(i)] $A(x).E^{u}(x)=E^{u}(T(x))$, $A(x).E^{s}(x)=E^{s}(T(x))$ for
a.e. $x\in M$;

\item[(ii)] there exist constants $C>0$ and $0<\tau <1$ such that 
\begin{eqnarray*}
\left\| A^{n}(x) \cdot v^{s}\right\| &\leq &C\tau ^{n}\left\| v^{s}\right\| , \\
\left\| A^{-n}(x) \cdot v^{u}\right\| &\leq &C\tau ^{n}\left\| v^{u}\right\|
\end{eqnarray*}
for a.e. $x\in X$, $v^{s}\in E^{s}(x)$, $v^{u}\in E^{u}(x)$ and $n\geq 0$;

\end{itemize}
\end{definition}

\begin{remark}
If $A$ is uniformly hyperbolic then
there is a constant $\alpha >0$ such that 
$\ang(E^{u}(x),E^{s}(x))\geq \alpha $ for a.e. $x\in X$;
see the proof of lemma~\ref{hiper bobo}.
\end{remark}

\begin{remark}
If $A$ is uniformly hyperbolic then it has
positive ($\geq -\log \tau $) Lyapunov exponent a.e. and the spaces 
$E^{u}$, $E^{s}$ in the definition coincide a.e. with the spaces given by
Oseledets theorem.
\end{remark}

\begin{remark}
The set $\mathcal{H}\subset L^{\infty }
(X,\sldoisr)$ of uniformly hyperbolic cocycles is open. This can be shown
by standard invariant cones techniques.
\end{remark}

We have the following result:

\begin{theorem B}
If $T$ is ergodic then there
is a residual set $\mathcal{R}\subset L^{\infty }(X,\sldoisr)$
such that for every $A\in \mathcal{R}$, either 
$A$ is uniformly hyperbolic or $\mathrm{LE}(A)=0$.
\end{theorem B}

From the above theorem we will deduce its continuous version. Now we suppose
that $X$ is a compact Hausdorff space, $\mu $ is a regular Borel probability
measure on $X$ and $T:(X,\mu )\hookleftarrow $ is an automorphism of $X.$ ($%
T $ is not assumed to be continuous). In this setting, we denote 
\[
C(X,\sldoisr)=\left\{ A:X\rightarrow \sldoisr%
\mbox{ continuous}\right\} ,
\]
endowed with the uniform convergence topology.
Then we have:

\begin{theorem C}
Let $X$, $\mu$ and $T$ be as above. If $T$ is ergodic then there is a residual set
$\mathcal{R}\subset C(X,\sldoisr)$ such that for every 
$A\in \mathcal{R}$, either $A$ is uniformly hyperbolic or 
$\mathrm{LE}(A)=0$.
\end{theorem C} 

%One says that a $T$-invariant set $Y\subset X$ is uniformly hyperbolic if
%the restricted cocycle $(Y,T|_{Y},A|_{Y})$ is uniformly hyperbolic.
%Theorem C has the following generalization to the non-ergodic case: 

%\begin{theorem D}
%Let $X$, $\mu$ and $T$ be as above. Then there is a residual set 
%$\mathcal{R}\subset C(X,\sldoisr)$
%such that for every $A\in \mathcal{R}$
%the following holds. If the invariant set 
%$\mathcal{O}^{+}=\{x\in X;\;\lambda ^{+}(A,x) > 0\}$
%has positive measure then there are invariant
%sets $H_{1}\subset H_{2}\subset \cdots \subset \mathcal{O}^{+}$
%covering $\mathcal{O}^{+}$ mod $0$ such that each set $H_{i}$
%is uniformly hyperbolic for $(T,A)$.
%\end{theorem D}

%An analogous generalization of theorem B to the non-ergodic case holds.

%%%%%%%%%%%%%%%%%%%%%%%%%%%%%%%%%%%%%%%%%%%%%%%%%%%%%%%%%%%%%%
\subsection{Proof of theorems B and C}
%%%%%%%%%%%%%%%%%%%%%%%%%%%%%%%%%%%%%%%%%%%%%%%%%%%%%%%%%%%%%%

\begin{proof}[Proof of theorem B]
The proof is similar to theorem A's, but it is easier in
several aspects.

Let $A\in L^{\infty }(X,\sldoisr)$ be a nonuniformly
hyperbolic cocycle with $\mathrm{LE}(A)>0$, $\delta >0$ and $\varepsilon >0$.
We have to show that there exists $\tilde{A}\in L^{\infty }(X,
\sldoisr)$ with $\| \tilde{A}-A\| _{\infty }<\varepsilon $
and $\mathrm{LE}(\tilde{A})<\delta$.

The first step is to prove an analogue of the Main Lemma.

\begin{lemma}\label{ML2}
Given $A$, $\delta$ and $\varepsilon$ as above, there exists a measurable
function $N:M\rightarrow \NN$ such that for a.e. $x\in X$ and every $%
n\geq N(x)$ there are $L_{0},\ldots ,L_{n-1}\in \sldoisr$
satisfying $\left\| L_{j}-A(T^{j}x)\right\| <\varepsilon $ and $\left\|
L_{n-1}\cdots L_{0}\right\| <e^{n\delta }$. Moreover, the matrices $L_{j}$
depend measurably on $x$ and $n$.
\end{lemma}

\begin{proof}
Since the proof of this lemma is essentially
contained in the proof of the Main Lemma, we will be brief. Let $\mathcal{O} \subset X$
be the full measure set given by Oseledets' theorem. For $m\in \NN$, let
$$
\Delta (x,m)=\frac{\left\| A^{m}(x)\mid _{E^{s}(x)}\right\| }{\left\|
A^{m}(x)\mid _{E^{u}(x)}\right\| }, \quad \text{for $x\in \mathcal{O}$},
$$
$$
\Gamma_m(A) = \left\{ x\in \mathcal{O}; \; \Delta(x,m) \geq 1/2 \right\}, \quad
\Omega_m(A) = \bigcup_{n\in\ZZ} T^n(\Gamma_m(A))
$$
Then, for every $m$, $\mu(\Omega_m(A))=1$ (same proof as lemma \ref{hiper bobo}).

Imitating the proof of lemma~\ref{Eu em Es}, one shows that
if $m$ is large enough then for every $y \in \Gamma_m(A)$ there exist
matrices $L'_{0}(y), \ldots, L'_{m-1}(y) \in \sldoisr$ such that
\begin{itemize}
\item[(i)] $\| L'_{j}(y) - A(T^{j} y) \| <\varepsilon$ and

\item[(ii)] $L'_{m-1}(y) \cdots L'_{0}(y) (E^{u}(y))=E^{s}(T^m y)$.
\end{itemize}
(These perturbations may be taken in the form
$L'_{j}(y) = A(T^{j}y)R_{\theta _{j}}$.) 

Proceeding as in the proof of lemma~\ref{ML simplificado}, one finds,
a measurable function $N: X \to \NN$ such that for a.e. $x\in X$ and
every $n \geq N(x)$, there is an integer $\ell \approx n/2$ such that
$y =T^\ell x \in \Gamma_m(A)$ and
$$
\| A^{n-\ell-m}(T^{\ell+m}(x)) \cdot L'_{m-1}(y) \cdots L'_0(y) \cdot A^{\ell}(x) \| <e^{n\delta}.
$$
The matrices $L_0, \ldots, L_{n-1}$ are defined in the obvious way:
$L_j = A(T^j(x))$ if $j<\ell$ or $j>\ell+m-1$;
and $L_j = L'_{j-\ell}(y)$ otherwise.
\end{proof}

We choose an integer $H$ and a set $B$ as in section 4.2. Let $\widehat{Q}$
be the Kakutani castle with base $B$. Since $T$ is ergodic and $\mu (B)>0$,
we have $\widehat{Q}=X$ mod 0. Again, let $Q$ be the castle consisting of all
the towers of $\widehat{Q}$ with heights of at most $3H$ floors. We have $\mu
(Q^{C})<3\delta ^{2}.$

Let $Q_{0}$ be the base of the castle $Q$. We apply lemma \ref{ML2} to a.e.
point in $Q_{0}$ to find $\tilde{A}\in L^{\infty }(X,\sldoisr )$ such that
\begin{itemize}
\item[(i)] $\| \tilde{A}-A\| _{\infty }<\varepsilon ,$

\item[(ii)] $\tilde{A}=A$ outside $Q,$

\item[(iii)] if $x\in Q_{0}$ and $n(x)$ is the height of the tower
containing $x$ then $\| \tilde{A}^{n}(x)\| <e^{n\delta }.$
\end{itemize}

Now let $N=\delta ^{-1}H$ and $G=\bigcap_{j=0}^{N-1}T^{-j}(Q)$. Since $Q^{C}$
is a castle with towers of height $\geq 3H,$ we have 
\begin{eqnarray*}
\mu (G^{C}) &=&\mu \left( \bigcup_{j=0}^{N-1}T^{-j}(Q^{C})\right) \leq \mu
(Q^{C})+N\mu \left( T^{-1}(Q^{C})-Q^{C}\right) \\
&\leq & \Big(1+\frac{N}{3H} \Big)\mu (Q^{C})<3\delta ^{2}+\delta <4\delta \, .
\end{eqnarray*}
Let $C>\left\| A\right\| _{\infty }+\varepsilon >1$. Then $\| 
\tilde{A}\| _{\infty }<C$ and if $x\in G$ then $\left\| \widetilde{%
A}^{N}(x)\right\| <C^{6H}e^{N\delta }$. Therefore
\begin{eqnarray*}
\mathrm{LE}(\tilde{A}) &\leq &\int_{G}\frac{1}{N}\log \left\| \widetilde{%
A}^{N}\right\| d\mu +\int_{G^{C}}\frac{1}{N}\log \| \tilde{A}%
^{N}\| d\mu  \\
&<&\left( 6\log C+1\right) \delta +(\log C)\mu (G^{C})<(10\log C+1)\delta 
\end{eqnarray*}
and theorem B is proved.
\end{proof}

\begin{remark}
An alternative proof of theorem B, without using
castles and the related estimates, can be given following Knill's \cite
{Knill} methods, using coboundary sets. This is done in~\cite{eu}.
\end{remark}

%In the proof of theorem B, ergodicity was hardly used.
%Indeed we have the stronger result:
%\begin{proposition}\label{nao ergodica}
%If $A \in L^\infty(X, \sldoisr)$ is such that every hyperbolic set has zero measure
%then for any $\delta >0$ and $\varepsilon >0$, there exists
%$\tilde{A}\in L^{\infty }(X,\sldoisr)$ with 
%$\| \tilde{A}-A\| _{\infty } < \varepsilon$ and $\mathrm{LE}(\tilde{A})<\delta$.
%\end{proposition}
%The only necessary modification in the proof is that instead of ``$\widehat{Q}=X$ mod 0''
%we would have ``$\widehat{Q}$ is a $\ T$-invariant set with 
%\mu (\widehat{Q})>1-\delta ^{2}$'' and then $\mu (G^{C})$ would be a bit larger.

\begin{proof}[Proof of theorem C]
It is possible to prove theorem C along the lines of the proof of theorem A,
but it is easier to deduce it from theorem B and Lusin's theorem.

Let $A\in C(X,\sldoisr)$ be a nonuniformly hyperbolic cocycle
with $\mathrm{LE}(A)>0,$ $\delta >0$ and $\varepsilon >0$. We have to show
that there exists $B\in C(X,\sldoisr)$ with $\left\|
B-A\right\| _{\infty }<\varepsilon $ and $\mathrm{LE}(B)<\delta .$

By theorem B, there exists $\tilde{A}\in L^{\infty }(X,
\sldoisr)$ near $A$ with $\mathrm{LE}(\tilde{A})=0$. Write $\tilde{A}%
=A.(I+J)$ with $J\in L^{\infty }(X,M(2,\RR ))$ close to $0$ in order that
if $J_{ij}(x)$ denote the entries of the matrix $J(x)$, then $\left\|
J_{ij}\right\| _{\infty }=\sup_{x}\left| J_{ij}(x)\right| <\varepsilon .$
Let $N\in \NN$ be such that $\frac{1}{N}\int_{X}\log \| \widetilde{A%
}^{N}\| d\mu <\delta $ and let $\gamma =N^{-1}\delta $.

By Lusin's theorem (see, for instance, \cite{Rudin}), there is $J^{\prime }$ 
$\in C(X,M(2,\RR ))$ with $\mu \{x;\;J^{\prime }\neq J\}<\gamma $ and $%
\left\| J_{ij}^{\prime }\right\| _{\infty }\leq \left\| J_{ij}\right\|
_{\infty }<\varepsilon $. Define $J_{ij}^{\prime \prime }=J_{ij}^{\prime }$
if $(i,j)\neq (1,1)$ and take $J_{11}^{\prime \prime }$ in order to have $%
\det (I+J^{\prime \prime })=1$. One can check that $\left\| J_{11}^{\prime
\prime }\right\| _{\infty }<2\varepsilon $ if we suppose 
$\varepsilon <1/3$. Thus
$\left\| J^{\prime \prime }\right\| _{\infty }\leq K\varepsilon$,
where $K$ is a constant. Let $B=A(I+J^{\prime \prime })$. Then $B\in C(X,%
\sldoisr)$ and 
\[
\left\| B-A\right\| _{\infty }\leq \left\| A\right\| _{\infty }\left\|
J^{\prime \prime }\right\| _{\infty }\leq K\left\| A\right\| _{\infty
}\varepsilon . 
\]
That is, $B$ is close to $A$.

Let $G_{0}=\{x;\;B(x)=\tilde{A}(x)\}$ and $G=%
\bigcap_{j=0}^{N-1}T^{-j}(G_{0})$. Then we have 
\[
\mu (G^{C})=\mu \left( \bigcup_{j=0}^{N-1}T^{-j}(G_{0}^{C})\right) <N\mu
(G_{0}^{C})<N\mu (\{x;\;J^{\prime }\neq J\})<\delta . 
\]
Besides, if $C>\left\| B\right\| _{\infty }$ then 
\begin{eqnarray*}
\mathrm{LE}(B) &\leq &\int_{G}\frac{1}{N}\log \left\| B^{N}\right\| d\mu
+\int_{G^{C}}\frac{1}{N}\log \left\| B^{N}\right\| d\mu \\
&\leq &\int_{G}\frac{1}{N}\log \| \tilde{A}^{N}\| d\mu +(\log
C)\mu (G^{C}) \\
&<&(1+\log C)\delta
\end{eqnarray*}
and we are done.
\end{proof}

%%%%%%%%%%%%%%%%%%%%%%%%%%%%%%%%%%%%%%%%%%%%%%%%%%%%%%%%%%
\subsection{The non-ergodic case}
%%%%%%%%%%%%%%%%%%%%%%%%%%%%%%%%%%%%%%%%%%%%%%%%%%%%%%%%%%

In the statements of theorems B and C, $T$ was assumed ergodic just for simplicity.
We will state without proof the generalization of these theorems to the non-ergodic case.

Again assume that $X$ is a compact Hausdorff space, $\mu$ is a regular Borel probability
measure on $X$ and $T:(X,\mu )\hookleftarrow $ is an automorphism of $X$.
One says that a $T$-invariant set $Y\subset X$ is uniformly hyperbolic if
the restricted cocycle $(Y,T|_{Y},A|_{Y})$ is uniformly hyperbolic.

\begin{theorem C'}
Let $X$, $\mu$ and $T$ be as above. Then there is a residual set 
$\mathcal{R}\subset C(X,\sldoisr)$
such that for every $A\in \mathcal{R}$
the following holds. If the invariant set 
$\mathcal{O}^{+}=\{x\in X;\;\lambda ^{+}(A,x) > 0\}$
has positive measure then there are invariant
sets $H_{1}\subset H_{2}\subset \cdots \subset \mathcal{O}^{+}$
covering $\mathcal{O}^{+}$ mod $0$ such that each set $H_{i}$
is uniformly hyperbolic for $(T,A)$.
\end{theorem C'}

The generalization of theorem B is entirely analogous.

%%%%%%%%%%%%%%%%%%%%%%%%%%%%%%%%%%%%%%%%%%%%%%%%%%%%%%%%%%
\subsection{Discontinuity of the Lyapunov exponents}
%%%%%%%%%%%%%%%%%%%%%%%%%%%%%%%%%%%%%%%%%%%%%%%%%%%%%%%%%%

Knill proves in \cite{Knill} that $\mathrm{LE}:L^{\infty }(X,
\sldoisr)\rightarrow \RR $ is discontinuous if $T$ is aperiodic. We will
analyze the continuous case.

Again, suppose that $X$ is a compact Hausdorff space and $\mu $ is a regular
Borel probability measure on $X$. Now let $T:(X,\mu )\hookleftarrow $ be an 
\emph{ergodic homeomorphism} of $X.$

By theorem C, $A\in C(X,\sldoisr)$ is a point of continuity of the function \linebreak
$\mathrm{LE}:C(X,\sldoisr)\rightarrow \RR $ if and only if
either $\mathrm{LE}(A)=0$ or $A$ is uniformly hyperbolic. Besides, it was
proved in \cite{Ruelle} that LE is even real-analytic when restricted to the
open set of uniformly hyperbolic cocycles.

Thus it is interesting to look for examples of continuous non-hyperbolic
cocycles with positive exponent for any given aperiodic ergodic system $%
(X,T,\mu )$. These examples are easily constructed if the system is not
uniquely ergodic, as is shown below.

\begin{proposition}
Let $(X,T,\mu )$ be as above. If $T$ is not uniquely ergodic then the
function $\mathrm{LE}:C(X,\sldoisr)\rightarrow \RR $ is
discontinuous.
\end{proposition}

\begin{proof}
Assume, without loss of generality, that the
support of $\mu $ is $X$. Take another invariant measure $\nu $. Take a
continuous function $h:X\rightarrow \RR $ such that $\int hd\mu \neq 0$
but $\int hd\nu =0$. Consider the diagonal cocycle 
\[
A(x)=\left( 
\begin{array}{ll}
e^{h(x)} & 0 \\ 
0 & e^{-h(x)}
\end{array}
\right) . 
\]
We have $\mathrm{LE}(A)=\left| \int hd\mu \right| $. Besides, for every $%
\varepsilon >0$ and $n_{0}>0$ there is $n>n_{0}$ such that the open set 
\[
\Bigg\{ x\in X;\;\frac{1}{n} \Big| \sum_{j=0}^{n-1}h(T^{j}x) \Big|
<\varepsilon \Bigg\} 
\]
is not empty and thus its $\mu$-measure is positive. This shows that $A$ is
not uniformly hyperbolic.
\end{proof}

For $T$ an irrational rotation of the circle some examples of continuous
nonuniform hyperbolic cocycles with positive exponent were given by Herman;
see \cite[\S 4]{Herman}. As a consequence, we obtain:

\begin{proposition}
For every irrational rotation $T$ of the $n$-torus $\mathbb{T}^{n}$, the function 
\[
\mathrm{LE}:C(\mathbb{T}^{n},\sldoisr)\rightarrow \RR  
\]
is discontinuous.
\end{proposition}

\begin{remark}
Clearly, it suffices to consider the case $n=1$.
\end{remark}

This proposition generalizes a previous result of Furman \cite{Furman},
which says that there is some irrational rotation such that $\mathrm{LE}$ is
discontinuous.

For a general aperiodic uniquely ergodic transformation $T$ it is an open question
whether such examples exist; see \cite{Furman}.

%%%%%%%%%%%%%%%%%%%%%%%%%%%%%%%%%%%%%%%%%%%%%%%%%%%%%%%%%
\medskip
\noindent \textit{Acknowledgments.}\,
%%%%%%%%%%%%%%%%%%%%%%%%%%%%%%%%%%%%%%%%%%%%%%%%%%%%%%%%%
I am grateful to my thesis advisor M.~Viana for numerous
conversations and suggestions and also for moral support.
Discussions with several other people at IMPA, specially
A.~Avila, were useful.
I would like to thank B.~Fayad for finding a gap in a preliminary 
version and F.~Abdenur, K.~Vixie and an anonymous referee for
corrections on the text.

%%%%%%%%%%%%%%%%%%%%%%%%%%%%%%%%%%%%%%%%%%%%%%%%%%%%%%%%%
%%%%%%%%%%%%%%%%%%%%%%%%%%%%%%%%%%%%%%%%%%%%%%%%%%%%%%%%%

\end{document}